\documentclass[12pt]{article}
\usepackage{amssymb}
\usepackage{amsmath}
\usepackage{latexsym}

\newtheorem{theorem}{Theorem}[section]
\newtheorem{proposition}[theorem]{Proposition}
\newtheorem{corollary}[theorem]{Corollary}
\newtheorem{lemma}[theorem]{Lemma}

\newtheorem{remark}[theorem]{Remark}

\begin{document}

\title{Free Bosonic Vertex Operator Algebras on Genus Two Riemann Surfaces~I}
\author{Geoffrey Mason\thanks{
Supported by the NSF, NSA, and the Committee on Research at the University
of California, Santa Cruz} \\
%EndAName
Department of Mathematics, \\
University of California Santa Cruz, \\
CA 95064, U.S.A. \and Michael P. Tuite\thanks{%
Partially supported by Millenium Fund, National University of Ireland, Galway%
} \\
%EndAName
School of Mathematics, Statistics and Applied Mathematics, \\
National University of Ireland, Galway \\
University Road,\\
Galway, Ireland.}
\maketitle

\begin{abstract}
We define the partition and $n$-point functions for a vertex operator
algebra on a genus two Riemann surface formed by sewing two tori together.
We obtain closed formulas for the genus two partition function for the
Heisenberg free bosonic string and for any pair of simple Heisenberg
modules. We prove that the partition function is holomorphic in the sewing
parameters on a given suitable domain and describe its modular properties
for the Heisenberg and lattice vertex operator algebras and a continuous
orbifolding of the rank two fermion vertex operator super algebra. We
compute the genus two Heisenberg vector $n$-point function and show that the
Virasoro vector one point function satisfies a genus two Ward identity for
these theories.
\end{abstract}

\tableofcontents

\normalfont

\section{Introduction}

\label{sect_intro} One of the most striking features of Vertex Operator
Algebras (VOAs) or chiral conformal field theory is the occurrence of \emph{%
elliptic functions} and \emph{modular forms}, manifested in the form of $n$%
-point correlation trace functions. This phenomenon has been present in
string theory since the earliest days e.g. \cite{GSW, P}. In mathematics it
dates from the Conway-Norton conjectures \cite{CN} proved by Borcherds (\cite%
{B1, B2}), and Zhu's important paper \cite{Z1}. Physically, we are dealing
with probability amplitudes corresponding to a complex torus (compact
Riemann surface of genus one) inflicted with $n$ punctures corresponding to
local fields (vertex operators). For a VOA $V=\oplus V_{n}$, the most
familiar correlation function is the $0$-point function, also called the 
\emph{partition function} or \emph{graded dimension} 
\begin{equation}
Z^{(1)}_{V}(q)=q^{-c/24}\sum_{n}\dim V_{n}q^{n},  \label{grdim}
\end{equation}%
($c$ is the central charge). An example which motivates much of the present
paper is that of a lattice theory $V_{L}$ associated to a positive-definite
even lattice $L$. Then $c$ is the rank of $L$ and 
\begin{equation}
Z_{V_{L}}^{(1)}(q)=\frac{\theta _{L}(q)}{\eta(q)^c},  \label{latticeg=1pfunc}
\end{equation}%
for Dedekind eta function $\eta(q)=q^{1/24}\prod_{n}(1-q^{n})$ and $\theta
_{L}(q)$ is the usual theta function of $L$. Both $\theta_{L}(q)$ and $%
\eta(q)^c$ are (holomorphic) elliptic modular forms of weight $c/2$ on a
certain congruence subgroup of $SL(2,\mathbb{Z})$, so that $Z_{V_{L}}$ is an
elliptic modular function of weight zero on the same subgroup. It is widely
expected that an analogous result holds for \emph{any rational} vertex
operator algebra, namely that $Z^{(1)}_{V}(q)$ is a modular function of
weight zero on a congruence subgroup of $SL(2,\mathbb{Z})$.

\medskip There are natural physical and mathematical reasons for wanting to
extend this picture to Riemann surfaces of \emph{higher genus}. In
particular, we want to know if there are natural analogs of (\ref{grdim})
and (\ref{latticeg=1pfunc}) for arbitrary rational vertex operator algebras
and arbitrary genus, in which \emph{genus g Siegel modular forms} occur.
This is considerably more challenging than the case of genus one. Many, but
not all, of the new difficulties that arise are already present at genus
two, and it is this case that we are concerned with in the present paper and
a companion paper \cite{MT4}. Our goal, then, is this: given a vertex
operator algebra $V$, to define the partition and $n$-point correlation
function on a compact Riemann surface of genus two which are associated to $%
V $, and study their convergence and automorphic properties. An overview of
aspects of this program is given in the Introduction to \cite{MT2}. Brief
discussions of some of our methods and results can also be found in \cite{T}, \cite{MT3} and 
\cite{MT6}.

\medskip The study of genus two (and higher) partition functions and
correlation functions has a long history in conformal field theory e.g. \cite%
{EO, FS, DP, So1, So2, BK, Kn, GSW, P} and, indeed, these ideas have heavily
influenced our approach. Likewise, in pure mathematics, other approaches
based on algebraic geometry have been been developed to describe $n$-point
correlation functions but \emph{not} the partition function e.g. \cite{TUY,
KNTY, Z2, U}. Our approach is constructively based \emph{only} on the
properties of a VOA in the spirit of Zhu's genus one theory \cite{Z1} with
no a priori assumptions made about the analytic or modular properties of
partition or $n$-point functions. Rather, in our approach, these genus two
objects are formally defined and are then proved to be analytic and modular
in appropriate domains for the VOAs considered.

\medskip In our approach, we define the genus two partition and $n$-point
functions in terms of genus one data coming from the VOA $V$. There are two
rather different ways to obtain a compact Riemann surface of genus two from
surfaces of genus one - one may sew two separate tori together, or self-sew
a torus (i.e. attach a handle). This is discussed at length in \cite%
{MT2} where we refer to these two schemes as the $\epsilon $- and $\rho $%
-formalism respectively. In the present paper we concentrate solely on
developing a theory of partition and $n$-point correlation functions in the $%
\epsilon$-formalism. We discuss the corresponding theory in the $\rho$%
-formalism in a companion paper \cite{MT4}.

\medskip The $\epsilon$-formalism developed in \cite{MT2} is reviewed
in Section~2 below. This is concerned with expressing a differential 2-form $%
\omega^{(2)}$ (the normalized differential of the second kind) in terms of a
pair of infinite matrices $A_i$, whose entries are quasi-modular forms
associated with the two sewn tori. This allows us to obtain explicit
expressions for genus two holomorphic one forms $\nu_1,\nu_2$ and the period
matrix $\Omega$ in terms of this genus one data. In particular, $\Omega$ is
determined by a holomorphic map 
\begin{equation}
\mathcal{D}^{\epsilon }\overset{F^{\epsilon }}{\longrightarrow }\mathbb{H}%
_{2},  \label{Fmaps1}
\end{equation}%
where for $g\geq 1,\mathbb{H}_{g}$ denotes the genus $g$ Siegel upper
half-space. Then $\mathcal{D}^{\epsilon }\subseteq \mathbb{H}_{1}\times 
\mathbb{H}_{1}\times \mathbb{C}$ is the domain consisting of triples $(\tau
_{1},\tau _{2},\epsilon )$ which correspond to a pair of complex tori of
modulus $\tau_{1},\tau _{2}$ sewn together by identifying two annular
regions via a sewing parameter $\epsilon$. This sewing produces a compact
Riemann surface of genus two, which assigns to each point of $\mathcal{D}%
^{\epsilon }$ the period matrix $\Omega $ of the sewn surface via the map $%
F^{\epsilon }$.

\medskip In Section~3 we introduce some graph-theoretic technology which
provides a convenient way of describing $\omega^{(2)},\ \nu_i$ and $\Omega$
in terms of the $\epsilon $-formalism. Similar graphical techniques are
employed later on as a means of computing the genus two partition function
and $n$-point functions for the free bosonic Heisenberg VOA and its modules.

\medskip Section~4 is a brief review of some necessary background on VOA
theory and the Li-Zamolodchikov or Li-Z metric. We assume throughout that
the Li-Z metric is unique and invertible (which follows if $V$ is simple 
\cite{Li}).

\medskip Section~5 develops a theory of $n$-point functions for VOAs on
Riemann surfaces of genus 0, 1 and 2 motivated by ideas in conformal field
theory. The Zhu theory \cite{Z1} of genus one $n$-point functions is
reformulated in this language in terms of the self-sewing of a Riemann
sphere to obtain a torus. We give a formal definition of genus two $n$-point
functions based on the given sewing formalism. We also emphasize the
interpretation of $n$-point functions in terms of formal differential forms.

The genus two partition function involves extending (\ref{Fmaps1}) to a
diagram 
\begin{equation*}
\begin{array}{ccc}
\mathcal{D}^{\epsilon } & \overset{F^{\epsilon }}{\longrightarrow } & \ 
\mathbb{H}_{2} \\ 
& \searrow & \ \downarrow \\ 
&  & \ \mathbb{C}%
\end{array}%
\end{equation*}%
where the partition function maps $\mathcal{D}^{\epsilon }\rightarrow 
\mathbb{C}$, and is defined purely in terms of genus one data coming from $V$%
. Explicitly, the genus two partition function of $V$ is \emph{a priori} a
formal power series in the variables $\epsilon ,q_{1},q_{2}$ (where as
usual, $q=e^{2\pi i\tau }$, etc.) given by 
\begin{equation}
Z_{V}^{(2)}(\tau _{1},\tau _{2},\epsilon ) =\sum_{n\geq 0}\epsilon
^{n}\sum_{u\in V_{[n]}}Z_{V}^{(1)}(u,\tau _{1})Z_{V}^{(1)}(\bar{u},\tau
_{2}),  \label{pfeps}
\end{equation}%
Here, $Z_{V}^{(1)}(u,\tau )$ is a genus one $1$-point function with $\bar{u}$
the Li-Z metric dual of $u$. The precise meaning of (\ref{pfeps}) together
with similar definitions for $n$-point functions, is given in Section~5.

\medskip In Sections~6 and 7 we investigate the case of the free bosonic
Heisenberg VOA $M$ and the expression corresponding to (\ref{pfeps}) for a
pair of simple $M$-modules. This later case is later used to analyze lattice
VOAs and the bosonized version of the rank two fermion Vertex Operator Super
Algebra. We find in all these cases that (\ref{pfeps}) is a holomorphic
function on $\mathcal{D}^{\epsilon }$. It is natural to expect that this
result holds in much wider generality. Section~6 is devoted to the
Heisenberg VOA $M$. In this case, holomorphy depends on an interesting new
formula for the genus two partition function. Namely, we prove (Theorem \ref%
{Theorem_Z2_boson}) by reinterpreting (\ref{pfeps}) in terms of certain
graphical expansion, that 
\begin{equation}
Z_{M}^{(2)}(\tau _{1},\tau _{2},\epsilon )=\frac{Z_{M}^{(1)}(\tau
_{1})Z_{M}^{(1)}(\tau _{2})}{\det (I-A_{1}A_{2})^{1/2}}.  \label{Z2M}
\end{equation}%
Here, the $A_{i}$ are the infinite matrices of Section~2 and $%
Z_{M}^{(1)}(\tau _{i})=1/\eta (q_{i})$. The infinite determinant that occurs
in (\ref{Z2M}) was introduced and discussed at length in \cite{MT2}. The
results obtained there are important here, as are the explicit computations
of genus one $1$-point functions obtained in \cite{MT1}. We also give in
Section~6 a product formula for the infinite determinant (Theorem \ref%
{Theorem_Z2_boson_prod}) which depends on the graphical interpretation of
the entries of the $A_{i}$.

The domain $\mathcal{D}^{\epsilon }$ admits the group $G_{0}=SL(2,\mathbb{Z}%
)\times SL(2,\mathbb{Z})$ as automorphisms (in fact, there is a larger
automorphism group $G$ that contains $G_{0}$ with index $2$). We show (cf.
Theorem \ref{Theorem_Z2_G}) that the partition function $Z_{M}^{(2)}(\tau
_{1},\tau _{2},\epsilon)$ is an automorphic form of weight $-1/2$ on $G$.
This is a bit imprecise in several ways: we have not explained here what the
automorphy factor is, and in fact this is an interesting point because it
depends on the map $F^{\epsilon }$. Similarly to the eta-function, there is
a $24$th root of unity, corresponding to a character of $G$, that intervenes
in the functional equation. These properties of $Z_{M}(\tau _{1},\tau
_{2},\epsilon )$ justify the idea that it should be thought of as the genus
two analog of $\eta (q)^{-1}$ in the $\epsilon $-formalism.

We conclude Section~6 by computing, by means of the graphical technique, the
genus two $n$-point function for $n$ Heisenberg vectors in terms of
symmetric tensor products of the differential 2-form $\omega^{(2)}$ in
Theorem \ref{theorem:Z2an}. This allows us to also find the Virasoro vector $%
1$-point function in terms of the genus two projective connection.

\medskip Section~7 is concerned with the genus two $n$-point function
associated with a pair of Heisenberg simple modules. We obtain a closed
formula for the partition function in Theorem \ref{Theorem_Z2_Malpha} and
the Heisenberg vector $n$-point function in terms of symmetric tensor
products of $\omega^{(2)}$ and $\nu_i$ in Theorem \ref{theorem:Z2anModule}.
We also derive a genus two Ward identity for the Virasoro vector 1-point
function in Proposition \ref{Prop:Z2omega alpha}. We apply these results in
Theorem \ref{Theorem_Z2_L_eps} to the case of a lattice VOA $V_L$ to find a
natural genus two generalization of (\ref{latticeg=1pfunc}), namely
\begin{equation}
\frac{Z_{V_L}^{(2)}(\tau _{1},\tau _{2},\epsilon )}{Z_{M}^{(2)}(\tau _{1},\tau _{2},\epsilon )}=\theta _{L}^{(2)}(\Omega),  \label{Z2Lattice}
\end{equation}%
where $\theta _{L}^{(2)}(\Omega )$ is the genus two Siegel theta function of the lattice $L$ . Similarly,
the Virasoro 1-point function obeys a Ward identity. Finally, we consider
the bosonized version of a continuous orbifolding of the rank two fermion
vertex super algebra to find the partition function is instead expressed in
terms of the genus two Riemann theta series.

\section{Genus Two Riemann Surface from Two Sewn Tori}

\label{sect_Riemann}

In this section we review some of the main results of \cite{MT2} relevant to
the present work. We review one of the two separate constructions of a genus
two Riemann surface discussed there based on a general sewing formalism due
to Yamada \cite{Y}. In this construction, which we refer to as the $\epsilon 
$-formalism, we parameterize a genus two Riemann surface by sewing together
two once-punctured tori. Then various genus two structures such as the
period matrix $\Omega$ can be determined in terms of genus one data. In
particular, $\Omega$ is described by an explicit formula which defines a
holomorphic map from a specified domain $\mathcal{D}^{\epsilon }$ into the
genus two Siegel upper half plane $\mathbb{H}_{2}$. This map is equivariant
under a suitable subgroup of $Sp(4,\mathbb{Z})$. We also review the
convergence and holomorphy of an infinite determinant that naturally arises
and which plays a dominant r\^{o}le later on.

\subsection{Some Elliptic Function Theory}

\label{subsect_Elliptic}

We begin with the definition of various modular and elliptic functions that
permeate this work \cite{MT1, MT2}. We define 
\begin{eqnarray}
P_{2}(\tau ,z) &=&\wp (\tau ,z)+E_{2}(\tau )  \notag \\
&=&\frac{1}{z^{2}}+\sum_{k=2}^{\infty }(k-1)E_{k}(\tau )z^{k-2},  \label{P2}
\end{eqnarray}%
where $\tau $ $\in \mathbb{H}_{1}$, the complex upper half-plane and where $%
\wp (\tau ,z)$ is the Weierstrass function and $E_{k}(\tau )$ is equal to $0$
for $k$ odd, and for $k$ even is the Eisenstein series 
\begin{equation}
E_{k}(\tau )=E_{k}(q)=-\frac{B_{k}}{k!}+\frac{2}{(k-1)!}\sum_{n\geq 1}\sigma
_{k-1}(n)q^{n}.  \label{Eisenk}
\end{equation}%
Here and below, we take $q=\exp (2\pi i\tau )$; $\sigma
_{k-1}(n)=\sum_{d\mid n}d^{k-1}$, and $B_{k}$ is a $k$th Bernoulli number
e.g. \cite{Se}. If $k\geq 4$ then $E_{k}(\tau )$ is a holomorphic modular
form of weight $k$ on $SL(2,\mathbb{Z})$ whereas $E_{2}(\tau )$ is a
quasi-modular form \cite{KZ, MT2}. We define $P_{1}(\tau ,z)$ by 
\begin{equation}
P_{1}(\tau ,z)=\frac{1}{z}-\sum_{k\geq 2}E_{k}(\tau )z^{k-1}.  \label{P1}
\end{equation}%
Noting $P_{2}=-\frac{d}{dz}P_{1}$ we define elliptic functions $P_{k}(\tau
,z)$ for $k\geq 3$ 
\begin{equation}
P_{k}(\tau ,z)=\frac{(-1)^{k-1}}{(k-1)!}\frac{d^{k-1}}{dz^{k-1}}P_{1}(\tau
,z).  \label{Pkdef}
\end{equation}

\noindent Define for $k,l\geq 1$ 
\begin{eqnarray}
C(k,l) &=&C(k,l,\tau )=(-1)^{k+1}\frac{(k+l-1)!}{(k-1)!(l-1)!}E_{k+l}(\tau ),
\label{Ckldef} \\
D(k,l,z) &=&D(k,l,\tau ,z)=(-1)^{k+1}\frac{(k+l-1)!}{(k-1)!(l-1)!}%
P_{k+l}(\tau ,z).  \label{Dkldef}
\end{eqnarray}%
The Dedekind eta-function is defined by

\begin{equation}
\eta (\tau )=q^{1/24}\prod_{n=1}^{\infty }(1-q^{n}).  \label{etafun}
\end{equation}

\subsection{The $\protect\epsilon $-Formalism for Sewing Two Tori}

\label{subsect_epsilon} Consider a compact Riemann surface $\mathcal{S}$ of
genus $2$ with canonical homology basis $a_{1},a_{2},b_{1},b_{2}$. There
exists two holomorphic 1-forms $\nu_{i}$, $i=1,2$ which we may normalize by 
\cite{FK} 
\begin{equation}
\oint_{a_{i}}\nu_{j}=2\pi i\delta_{ij}.  \label{norm}
\end{equation}%
These forms can also be defined via the unique singular bilinear two form $%
\omega^{(2)}$, known as the \emph{normalized differential of the second kind}%
. It is defined by the following properties \cite{FK, Y}: 
\begin{equation}
\omega^{(2)}(x,y)=(\frac{1}{(x-y)^{2}}+\text{regular terms})dxdy
\label{omegag}
\end{equation}%
for any local coordinates $x,y$, with normalization 
\begin{equation}
\int_{a_{i}}\omega^{(2)}(x,\cdot )=0,  \label{nugnorm}
\end{equation}%
for $i=1,2$. Using the Riemann bilinear relations, one finds that 
\begin{equation}
\nu_{i}(x)=\oint_{b_{i}}\omega^{(2)}(x,\cdot ),  \label{nui}
\end{equation}%
with $\nu_{i}$ normalized as in (\ref{norm}). The genus $2$ period matrix $%
\Omega$ is then defined by 
\begin{equation}
\Omega_{ij}=\frac{1}{2\pi i}\oint_{b_{i}}\nu_{j}\quad  \label{period}
\end{equation}%
for $i,j=1,2$. One further finds that $\Omega \in \mathbb{H}_{2}$, the
Siegel upper half plane.

We now review a general method due to Yamada \cite{Y} and discussed at
length in \cite{MT2} for calculating $\omega^{(2)}(x,y)$, $\nu_i(x)$ and $%
\Omega_{ij} $ on the genus two Riemann surface formed by sewing together two
tori $\mathcal{S}_{a} $ for $a=1,2$. We shall sometimes refer to $\mathcal{S}%
_{1}$ and $\mathcal{S}_{2}$ as the left and right torus respectively.
Consider an oriented torus $\mathcal{S}_{a}=\mathbb{C}/\Lambda _{a}$ with
lattice $\Lambda _{a}=2\pi i(\mathbb{Z}\tau _{a}\oplus \mathbb{Z})$ for $%
\tau _{a}\in \mathbb{H}_{1}$. For local coordinate $z_{a}\in \mathbb{C}%
/\Lambda _{a}$ consider the closed disk $\left\vert z_{a}\right\vert \leq
r_{a}$ which is contained in $\mathcal{S}_{a}$ provided\ $r_{a}<\frac{1}{2}%
D(q_{a})$ where 
\begin{equation*}
D(q_{a})=\min_{\lambda \in \Lambda _{a},\lambda \neq 0}|\lambda |,
\end{equation*}%
is the minimal lattice distance. Introduce a complex sewing parameter $%
\epsilon $ where $|\epsilon |\leq r_{1}r_{2}<\frac{1}{4}D(q_{1})D(q_{2})$
and excise the disk $\{z_{a},\left\vert z_{a}\right\vert \leq |\epsilon |/r_{%
\bar{a}}\}$ centered at $z_{a}=0$ to form a punctured torus 
\begin{equation*}
\hat{\mathcal{S}}_{a}=\mathcal{S}_{a}\backslash \{z_{a},\left\vert
z_{a}\right\vert \leq |\epsilon |/r_{\bar{a}}\},
\end{equation*}%
where we use the convention 
\begin{equation}
\overline{1}=2,\quad \overline{2}=1.  \label{abar}
\end{equation}%
Defining the annulus

\begin{equation}
\mathcal{A}_{a}=\{z_{a},|\epsilon |/r_{\bar{a}}\leq \left\vert
z_{a}\right\vert \leq r_{a}\}\subset \hat{\mathcal{S}}_{a},
\label{annuli_eps}
\end{equation}%
we identify $\mathcal{A}_{1}$ with $\mathcal{A}_{2}$ via the sewing relation 
\begin{equation}
z_{1}z_{2}=\epsilon.  \label{pinch}
\end{equation}
The genus two Riemann surface is parameterized by the domain 
\begin{equation}
\mathcal{D}^{\epsilon }=\{(\tau _{1},\tau _{2},\epsilon )\in \mathbb{H}_{1}%
\mathbb{\times H}_{1}\mathbb{\times C}\ |\ |\epsilon |<\frac{1}{4}%
D(q_{1})D(q_{2})\}.  \label{Deps}
\end{equation}

\begin{center}
\begin{picture}(300,100)

%left surface
\put(50,50){\qbezier(10,-20)(-30,0)(10,20)}% left left
\put(50,52){\qbezier(10,18)(50,35)(90,18)}% left upper

\put(50,48){\qbezier(10,-18)(50,-35)(90,-18)}%left lower

\put(45,50){\qbezier(25,0)(45,17)(60,0)}%upper
\put(45,50){\qbezier(20,2)(45,-17)(65,2)}%lower

%right surface

\put(175,52){\qbezier(10,18)(50,35)(90,18)}%right upper
\put(175,50){\qbezier(90,20)(130,0)(90,-20)}%right right
\put(175,48){\qbezier(10,-18)(50,-35)(90,-18)}% right lower

\put(200,50){\qbezier(25,0)(45,17)(60,0)}%upper
\put(200,50){\qbezier(20,2)(45,-17)(65,2)}%lower
% left annulus centered at (140,50)

\put(140,50){\circle{16}}
\put(140,50){\circle{40}}

% z_1=0label
\put(140,50){\vector(-1,-2){0}}%arrow
\put(50,50){\qbezier(90,0)(100,15)(90,30)}%
\put(140,90){\makebox(0,0){$z_1=0$}}

% line and r1 label
\put(140,50){\line(-1,1){14.1}}
\put(127,55){\makebox(0,0){$r_1$}}

% line and eps/r2 label
\put(140,50){\line(1,0){8}}
\put(145,50){\vector(1,4){0}}%arrow
\put(55,20){\qbezier(90,4)(85,17)(90,30)}%
\put(150,15){\makebox(0,0){$|\epsilon|/r_2$}}

%Sg1 label
\put(25,50){\makebox(0,0){$\mathcal{S}_1$}}

%right annulus centred at (185,50)

\put(185,50){\circle{16}}
\put(185,50){\circle{40}}

% z_2=0label
\put(185,50){\vector(1,-2){0}}%arrow
\put(95,50){\qbezier(90,0)(80,15)(90,30)}%
\put(185,90){\makebox(0,0){$z_2=0$}}

% line and r2 label
\put(185,50){\line(-1,-1){14.1}}
\put(173,46){\makebox(0,0){$r_2$}}

% line and eps/r1 label
\put(185,50){\line(1,0){8}}
\put(190,50){\vector(1,4){0}}%arrow
\put(90,20){\qbezier(100,4)(95,17)(100,30)}%
\put(190,15){\makebox(0,0){$ |\epsilon|/r_1$}}

%Sg2 label
\put(300,50){\makebox(0,0){$\mathcal{S}_2$}}

\end{picture}

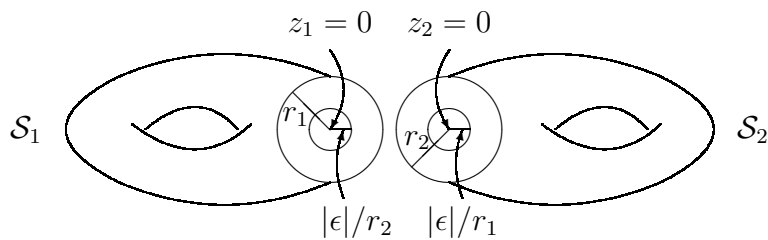
{\small Fig.~1 Sewing Two Tori}
\end{center}

We next introduce the infinite dimensional matrix $A_{a}(\tau _{a},\epsilon
)=(A_{a}(k,l,\tau_{a},\epsilon ))$ for $k,l\geq 1$ where 
\begin{equation}
A_{a}(k,l,\tau _{a},\epsilon )=\frac{\epsilon ^{(k+l)/2}}{\sqrt{kl}}%
C(k,l,\tau _{a}).  \label{Akldef}
\end{equation}%
The matrices $A_1,A_2$ play a dominant role both here and in our later
discussion of the free bosonic VOA and its modules on a genus two Riemann
surface. In particular, the matrix $I-A_{1}A_{2}$ and $\det (I-A_{1}A_{2})$
(where $I$ denotes the infinite identity matrix) play an important role
where $\det (I-A_{1}A_{2})$ is defined by 
\begin{eqnarray}
\log \det (I-A_{1}A_{2}) &=&\mathrm{Tr}\log (I-A_{1}A_{2})  \notag \\
&=&-\sum_{n\geq 1}\frac{1}{n}\mathrm{Tr}((A_{1}A_{2})^{n}).  \label{logdet}
\end{eqnarray}%
One finds

\begin{theorem}
\label{Theorem_A1A2} \ \ 

(a) (op. cite., Proposition 1) The infinite matrix 
\begin{equation}
(I-A_{1}A_{2})^{-1}=\sum_{n\geq 0}(A_{1}A_{2})^{n},  \label{I_minus_A1A2}
\end{equation}%
is convergent for $(\tau _{1},\tau _{2},\epsilon )\in \mathcal{D}^{\epsilon}$.

(b) (op. cite., Theorem 2 \& Proposition 3) $\det (I-A_{1}A_{2})$ is
non-vanishing and holomorphic for $(\tau _{1},\tau _{2},\epsilon )\in 
\mathcal{D}^{\epsilon }$. $\square $
\end{theorem}

The bilinear two form $\omega^{(2)}(x,y)$, the holomorphic one forms $%
\nu_i(x)$ and the period matrix $\Omega_{ij}$ are given in terms of the
matrices $A_a$ and holomorphic one forms on the punctured torus $\hat{%
\mathcal{S}}_{a}$ given by 
\begin{equation}
a_{a}(k,x)=\sqrt{k}\epsilon ^{k/2}P_{k+1}(\tau _{a},x)dx.  \label{eq:1forma}
\end{equation}
Letting $a_{a}(x)$, $a_{a}^T(x)$ denote the infinite row, respectively
column vector with elements (\ref{eq:1forma}) we have:

\begin{theorem}
\label{Theorem_om_eps}(op. cite., Lemma 2, Proposition 1, Theorem 4) 
\begin{equation}
\omega^{(2)}(x,y)=\left\{ 
\begin{array}{ll}
P_2(\tau_a,x-y)dxdy
+ a_{a}(x)A_{\bar{a}}(I-A_a A_{\bar{a}})^{-1}a_{a}^{T}(y),
& x, y\in \hat{\mathcal{S}}_{a}, \\ 
-a_{a}(x)(I-A_{\bar{a}}A_a)^{-1}a_{\bar{a}}^{T}(y), & x\in \hat{\mathcal{S}}%
_{a},\ y\in \hat{\mathcal{S}}_{\bar{a}}.%
\end{array}
\right.  \label{om_eps}
\end{equation}
$\square $
\end{theorem}

Applying (\ref{nui}) we then find (op. cite., Theorem 4) 
\begin{equation}
\nu_a(x)=\left\{ 
\begin{array}{ll}
dx+\epsilon^{1/2}(a_{a}(x)A_{\bar{a}}(I-A_a A_{\bar{a}})^{-1})(1) & x\in 
\hat{\mathcal{S}}_{a}, \\ 
-\epsilon^{1/2}(a_{\bar{a}}(x)(I-A_a A_{\bar{a}})^{-1})(1) & x\in \hat{%
\mathcal{S}}_{\bar{a}},%
\end{array}
\right.  \label{nui_eps}
\end{equation}
where $(1)$ refers to the $(1)$-entry of a vector. Furthermore applying (\ref%
{period}) we have

\begin{theorem}
\label{Theorem_period_eps}(op. cite., Theorem 4) The $\epsilon$-formalism
determines a holomorphic map 
\begin{eqnarray}
F^{\epsilon }:\mathcal{D}^{\epsilon } &\rightarrow &\mathbb{H}_{2},  \notag
\\
(\tau _{1},\tau _{2},\epsilon ) &\mapsto &\Omega (\tau _{1},\tau
_{2},\epsilon ),  \label{Fepsmap}
\end{eqnarray}%
where $\Omega =\Omega (\tau _{1},\tau _{2},\epsilon )$ is given by 
\begin{eqnarray}
2\pi i\Omega _{11} &=&2\pi i\tau _{1}+\epsilon
(A_{2}(I-A_{1}A_{2})^{-1})(1,1),  \label{Om11eps} \\
2\pi i\Omega _{22} &=&2\pi i\tau _{2}+\epsilon
(A_{1}(I-A_{2}A_{1})^{-1})(1,1),  \label{Om22eps} \\
2\pi i\Omega _{12} &=&-\epsilon (I-A_{1}A_{2})^{-1}(1,1).  \label{Om12eps}
\end{eqnarray}%
Here $(1,1)$ refers to the $(1,1)$-entry of a matrix. \hfill $\square $
\end{theorem}

$\mathcal{D}^{\epsilon }$ is preserved under the action of $G\simeq (SL(2,%
\mathbb{Z})$ $\times SL(2,\mathbb{Z}))\rtimes \mathbb{Z}_{2}$, the direct
product of two copies of $SL(2,\mathbb{Z})$ (the left and right torus
modular groups) which are interchanged upon conjugation by an involution $%
\beta $ as follows%
\begin{eqnarray}
\gamma _{1}.(\tau _{1},\tau _{2},\epsilon ) &=&(\frac{a_{1}\tau _{1}+b_{1}}{%
c_{1}\tau _{1}+d_{1}},\tau _{2},\frac{\epsilon }{c_{1}\tau _{1}+d_{1}}), 
\notag \\
\gamma _{2}.(\tau _{1},\tau _{2},\epsilon ) &=&(\tau _{1},\frac{a_{2}\tau
_{2}+b_{2}}{c_{2}\tau _{2}+d_{2}},\frac{\epsilon }{c_{2}\tau _{2}+d_{2}}), 
\notag \\
\beta .(\tau _{1},\tau _{2},\epsilon ) &=&(\tau _{2},\tau _{1},\epsilon ),
\label{GDeps}
\end{eqnarray}%
for $(\gamma _{1},\gamma _{2})\in SL(2,\mathbb{Z})\times SL(2,\mathbb{Z})$
with $\gamma _{i}=\left( 
\begin{array}{cc}
a_{i} & b_{i} \\ 
c_{i} & d_{i}%
\end{array}%
\right) $. There is a natural injection $G\rightarrow Sp(4,\mathbb{Z})$ in
which the two $SL(2,\mathbb{Z})$ subgroups are mapped to 
\begin{equation}
\Gamma _{1}=\left\{ \left[ 
\begin{array}{cccc}
a_{1} & 0 & b_{1} & 0 \\ 
0 & 1 & 0 & 0 \\ 
c_{1} & 0 & d_{1} & 0 \\ 
0 & 0 & 0 & 1%
\end{array}%
\right] \right\} ,\;\Gamma _{2}=\left\{ \left[ 
\begin{array}{cccc}
1 & 0 & 0 & 0 \\ 
0 & a_{2} & 0 & b_{2} \\ 
0 & 0 & 1 & 0 \\ 
0 & c_{2} & 0 & d_{2}%
\end{array}%
\right] \right\} ,  \label{Gamma1Gamma2}
\end{equation}%
and the involution is mapped to 
\begin{equation}
\beta =\left[ 
\begin{array}{cccc}
0 & 1 & 0 & 0 \\ 
1 & 0 & 0 & 0 \\ 
0 & 0 & 0 & 1 \\ 
0 & 0 & 1 & 0%
\end{array}%
\right] .  \label{betagen}
\end{equation}%
Thus as a subgroup of $Sp(4,\mathbb{Z})$, $G$ also has a natural action on
the Siegel upper half plane $\mathbb{H}_{2}$ where for $\gamma =\left( 
\begin{array}{ll}
A & B \\ 
C & D%
\end{array}%
\right) \in Sp(4,\mathbb{Z})$ 
\begin{equation}
\gamma .\Omega {=(A\Omega +B)(C\Omega +D)^{-1},}  \label{eq: modtrans}
\end{equation}%
One then finds

\begin{theorem}
\label{TheoremGequiv} (op. cit., Theorem 5) $F^{\epsilon }$ is equivariant
with respect to the action of $G$ i.e. there is a commutative diagram for $%
\gamma \in G$, 
\begin{equation*}
\begin{array}{ccc}
\mathcal{D}^{\epsilon } & \overset{F^{\epsilon }}{\rightarrow } & \mathbb{H}%
_{2} \\ 
\gamma \downarrow &  & \downarrow \gamma \\ 
\mathcal{D}^{\epsilon } & \overset{F^{\epsilon }}{\rightarrow } & \mathbb{H}%
_{2}%
\end{array}%
\end{equation*}
$\square$
\end{theorem}

\section{Graphical expansions}

\label{sect_graph}

\subsection{Rotationless and Chequered Cycles}

\label{subsect_cycles}

We set up some notation and discuss certain types of labeled graphs. These
arise directly from consideration of the terms that appear in the
expressions for $\omega^{(2)}(x,y)$, $\nu_i(x)$ and $\Omega_{ij}$ reviewed
in the last Section, and will later play an important r\^{o}le in the
analysis of genus two partition functions for vertex operator algebras.

\bigskip Next we introduce the notion of a \emph{chequered cycle} as a
(clockwise) oriented, labeled polygon $L$ with $2n$ nodes for some integer $%
n\geq 0$, and nodes labeled by arbitrary positive integers. Moreover, edges
carry a label $1$ or $2$ which alternate as one moves around the polygon.

\bigskip

\begin{center}
\begin{picture}(250,80)

%\graphpaper(0,0)(250,100)

\put(100,50){\line(1,2){10}}
\put(90,50){\makebox(0,0){$i_1$}}
\put(100,50){\circle*{4}}
\put(102,63){\makebox(0,0){\scriptsize 1}}
\put(107,63){\vector(1,2){0}}%arrow

\put(110,70){\line(1,0){20}}
\put(105,78){\makebox(0,0){$i_2$}}
\put(110,70){\circle*{4}}
\put(120,75){\makebox(0,0){\scriptsize 2}}
\put(122,70){\vector(1,0){0}}%arrow

\put(130,70){\line(1,-2){10}}
\put(139,78){\makebox(0,0){$i_3$}}
\put(130,70){\circle*{4}}
\put(139,63){\makebox(0,0){\scriptsize 1}}
\put(137,57){\vector(1,-2){0}}%arrow

\put(140,50){\line(-1,-2){10}}
\put(150,50){\makebox(0,0){$i_4$}}
\put(140,50){\circle*{4}}
\put(138,36){\makebox(0,0){\scriptsize 2}}
\put(134,38){\vector(-1,-2){0}}%arrow

\put(110,30){\line(1,0){20}}
\put(135,20){\makebox(0,0){$i_5$}}
\put(130,30){\circle*{4}}
\put(120,25){\makebox(0,0){\scriptsize 1}}
\put(117,30){\vector(-1,0){0}}%arrow

\put(100,50){\line(1,-2){10}}
\put(105,20){\makebox(0,0){$i_6$}}
\put(110,30){\circle*{4}}
\put(102,36){\makebox(0,0){\scriptsize 2}}
\put(104,42){\vector(-1,2){0}}%arrow
\end{picture}

{\small Fig.~2 Chequered Cycle}
\end{center}

We call a node with label $1$ \emph{distinguished} if its abutting edges are
of type $\overset{2}{\longrightarrow }\overset{1}{\bullet }\overset{1}{%
\longrightarrow }$. Set 
\begin{eqnarray}  \label{mathcalRdef}
\mathcal{R} &=& \{%
\mbox{isomorphism classes of rotationless chequered cycles
} \}, \\
\mathcal{R}_{21} &=&\{ 
\mbox{isomorphism classes of rotationless  chequered
cycles }  \notag \\
&&\mbox{ with a  distinguished node}\},  \notag \\
\mathcal{L}_{21} &=&\{\mbox{isomorphism classes of chequered cycles  with a}
\notag \\
&&\mbox{ \emph{unique} distinguished node}\},  \notag
\end{eqnarray}

\medskip

Let $S$ be a commutative ring and $S[t]$ the polynomial ring with
coefficients in $S$. Let $M_{1}$ and $M_{2}$ be infinite matrices with $%
(k,l) $-entries 
\begin{equation}
M_{a}(k,l)=t^{k+l}s_{a}(k,l)  \label{eq: abstractAmatrix}
\end{equation}%
for $a=1,2$ and $k,l\geq 1$, where $s_{a}(k,l)\in S$. Given this data, we
define a map, or \emph{weight function}, 
\begin{equation*}
\zeta :\{\mbox{chequered cycles}\}\longrightarrow S[t]
\end{equation*}%
as follows: if $L$ is a chequered cycle then $L$ has edges $E$ labeled as \ $%
\overset{k}{\bullet }\overset{a}{\longrightarrow }\overset{l}{\bullet }$.
Then set $\zeta (E)=M_{a}(k,l)$ and 
\begin{equation}
\zeta (L)=\prod \zeta (E)  \label{eq: eqomega}
\end{equation}%
where the product is taken over all edges of $L$.

\bigskip It is useful to also introduce a variation on the theme of
chequered polygons, namely \emph{\ oriented chequered necklaces}. These are
connected graphs with $n\geq 3$ nodes, $(n-2)$ of which have valency $2$ and
two of which have valency $1$ (these latter are the \emph{end nodes})
together with an orientation, say from left to right. There is also a
degenerate necklace $N_{0}$ with a single node and no edges. As before,
nodes are labeled with arbitrary positive integers and edges are labeled
with an index $1$ or $2$ which alternate along the necklace. For such a
necklace $N$, we define the weight function $\zeta (N)$ as a product of edge
weights as in (\ref{eq: eqomega}), with $\zeta (N_{0})=1$.

\medskip Among all chequered necklaces there is a distinguished set for
which both end nodes are labeled by $1$. There are four types of such
chequered necklaces, which may be further distinguished by the labels of the
two edges at the extreme left and right. Using the convention (\ref{abar})
we say that the chequered necklace 
\begin{eqnarray*}
\overset{1}{\bullet }\overset{a}{\longrightarrow }\overset{i}{\bullet }%
&\ldots& \overset{j}{\bullet }\overset{b}{\longrightarrow }\overset{1}{\bullet 
}\\
&\mbox{\small {Fig.~3}}&
\end{eqnarray*}

\noindent is of \emph{type $ab$} for $a,b\in \{1,2\}$, and set 
\begin{eqnarray}  \label{Nabdef}
\mathcal{N}_{ab} &=&\{\mbox{isomorphism classes of oriented chequered} 
\notag \\
&&\mbox{ \ necklaces of type}\ ab\},
\end{eqnarray}
\begin{eqnarray}  \label{omegaabdef}
\zeta _{ab} &=&\sum_{N\in \mathcal{N}_{ab}}\zeta (N).
\end{eqnarray}

\subsection{Necklace Graphical Expansions for $\protect\omega^{(2)}$, $%
\protect\nu_i$ and $\Omega_{ij}$}

\label{subsect_Omega_graph}

We now apply the formalism of the previous Subsection to the expressions for 
$\omega^{(2)}(x,y)$, $\nu_i(x)$ and $\Omega_{ij}$ in the $\epsilon $%
-formalism reviewed in Section 2. We begin with the period matrix $%
\Omega_{ij}$. Here the ring $S$ is taken to be the product $S_{1}\times
S_{2} $ where for $a=1,2$, $S_{a}$ is the ring of quasi-modular forms $%
\mathbb{C}[E_{2}(\tau _{a}), E_{4}(\tau_{a}),E_{6}(\tau _{a})]$, and $t =
\epsilon^{1/2}$. The matrices $M_a$ are taken to be the $A_a$ defined in (%
\ref{Akldef}). Thus 
\begin{equation}
s_{a}(k,l)=\frac{C(k,l,\tau _{a})}{\sqrt{kl}},  \label{eq: subst}
\end{equation}%
and for the edge $E$ labeled as $\overset{k}{\bullet}\overset{a}{%
\longrightarrow}\overset{l}{\bullet}$ we have 
\begin{eqnarray}  \label{omegaval}
\zeta (E)=A_{a}(k,l).
\end{eqnarray}

\medskip Recalling the notation (\ref{omegaabdef}), we find

\begin{proposition}
\label{Propepsperiodgraph} (\cite{MT2}, Proposition 4) For $a=1,2$ 
\begin{eqnarray*}
\Omega _{aa} &=&\tau _{a}+\frac{\epsilon }{2\pi i}\zeta _{\bar{a}\bar{a}}, \\
\Omega _{a\bar{a}} &=&-\frac{\epsilon }{2\pi i}\zeta _{\bar{a}a}.\ \ \hfill
\square
\end{eqnarray*}
\end{proposition}

Furthermore, in the notation of Section \ref{subsect_cycles} we have

\begin{proposition}
\label{Prop_Om12_R21expansion} 
\begin{equation}
\zeta_{12}=\zeta_{21}=\prod_{L\in \mathcal{R}_{21}}(1-\zeta(L))^{-1}.
\label{eq: omega_12prod}
\end{equation}
\end{proposition}

\noindent Beyond the intrinsic interest of this product formula, our main
use of it will be to provide an alternate proof of Theorem \ref{Theorem_Z2_G}
below. We therefore relegate the proof to Proposition \ref{Prop_R21expansion}
to the Appendix.

We can similarly obtain necklace graphical expansions for the bilinear form $%
\omega^{(2)}(x,y)$ and the holomorphic one forms $\nu_i(x)$. We introduce
further distinguished valence one nodes labeled by $1,x$ for $x\in \hat{%
\mathcal{S}_a}$. The set of edges $\{E\}$ is augmented by edges with weights
defined by: 
\begin{eqnarray}
\zeta(\overset{1,x}{\bullet}\overset{a}{\longrightarrow}\overset{1,y}{\bullet%
}) &=&P_2(\tau_a,x-y),\ x,y\in \hat{\mathcal{S}_a},  \notag \\
\zeta(\overset{1,x}{\bullet}\overset{a}{\longrightarrow}\overset{k}{\bullet}%
) =\zeta(\overset{k}{\bullet}\overset{a}{\longrightarrow}\overset{1,x}{%
\bullet}) &=&\sqrt{k}\epsilon^{k/2}P_{k+1}(\tau_a,x),\ x\in \hat{\mathcal{S}%
_a},  \label{eq: zeta1a2}
\end{eqnarray}
for elliptic functions (\ref{Pkdef}).

Similarly to (\ref{Nabdef}) we consider chequered necklaces where one or
both end points are $1,x$-type labeled nodes. We thus define for $x\in \hat{%
\mathcal{S}}_{ a}$ and $y\in \hat{\mathcal{S}}_{b}$ three isomorphism
classes of oriented chequered necklaces denoted $\mathcal{N}_{ab}^{x,1}$, $%
\mathcal{N}_{ab}^{1,y}$ and $\mathcal{N}_{ab}^{x,y}$ with the following
respective configurations 
\begin{eqnarray}
\{\overset{1,x}{\bullet }\overset{a}{\longrightarrow }\overset{i}{\bullet }%
\ldots \overset{j}{\bullet}\overset{b}{\longrightarrow }\overset{1}{\bullet }%
\},  \label{Nabx1def} \\
\{\overset{1}{\bullet }\overset{a}{\longrightarrow }\overset{i}{\bullet }%
\ldots \overset{j}{\bullet}\overset{b}{\longrightarrow }\overset{1,y}{%
\bullet }\},  \label{Nab1ydef} \\
\{\overset{1,x}{\bullet }\overset{a}{\longrightarrow }\overset{i}{\bullet }%
\ldots \overset{j}{\bullet}\overset{b}{\longrightarrow }\overset{1,y}{%
\bullet }\}.  \label{Nabxydef}
\end{eqnarray}
Let $\zeta _{ab}^{x,1}$, $\zeta _{ab}^{1,y}$ and $\zeta _{ab}^{x,y}$ denote
the respective sum of the weights for each class. Comparing to (\ref{om_eps}%
) and (\ref{nui_eps}) and applying (\ref{nui}) we find the following
graphical expansions for the bilinear form $\omega^{(2)}(x,y)$ and the
holomorphic one forms $\nu_i(x)$

\begin{proposition}
\label{Propomeganugraph} For $a=1,2$ 
\begin{eqnarray}
\omega^{(2)}(x,y)&=&\left\{ 
\begin{array}{ll}
\zeta _{aa}^{x,y}dxdy & x, y\in \hat{\mathcal{S}}_{a}, \\ 
-\zeta _{a\bar{a}}^{x,y}dxdy & x\in \hat{\mathcal{S}}_{a},\ y\in \hat{%
\mathcal{S}}_{\bar{a}},%
\end{array}
\right.  \label{om_graph} \\
\nu_a(x)&=&\left\{ 
\begin{array}{ll}
(1+\epsilon^{1/2} \zeta _{aa}^{x,1})dx & x\in \hat{\mathcal{S}}_{a}, \\ 
-\epsilon^{1/2}\zeta _{\bar{a}a}^{x,1}dx & x\in \hat{\mathcal{S}}_{\bar{a}}.%
\end{array}
\right.  \label{nui_graph}
\end{eqnarray}
\end{proposition}

\section{Vertex Operator Algebras and the Li-Zamolodchikov Metric}

\label{sect_VOALIZ}

\subsection{Vertex Operator Algebras}

\label{subsect_VOA}

We review some relevant aspects of vertex operator algebras (\cite{FHL, FLM,
Ka, LL, MN}). A vertex operator algebra (VOA) is a quadruple $(V,Y,\mathbf{1}%
,\omega )$ consisting of a $\mathbb{Z}$-graded complex vector space $%
V=\bigoplus_{n\in \mathbb{Z}}V_{n}$, a linear map $Y:V\rightarrow (\mathrm{%
End}\ V)[[z,z^{-1}]]$, for formal parameter $z$, and a pair of distinguished
vectors (states): the vacuum $\mathbf{1}\in V_{0}$ , and the conformal
vector $\omega \in V_{2}$. For each state $v\in V$ the image under the $Y$
map is the vertex operator

\begin{equation}
Y(v,z)=\sum_{n\in \mathbb{Z}}v(n)z^{-n-1},  \label{Ydefn}
\end{equation}%
with modes $v(n)\in \mathrm{End}V$ where $\mathrm{Res}_{z=0}z^{-1}Y(v,z)%
\mathbf{1}=v(-1)\mathbf{1}=v$. Vertex operators satisfy the Jacobi identity
or equivalently, operator locality or Borcherds's identity for the modes
(loc. cit.).

The vertex operator for the conformal vector $\omega $ is defined as

\begin{equation*}
Y(w,z)=\sum_{n\in \mathbb{Z}}L(n)z^{-n-2}.
\end{equation*}%
The modes $L(n)$ satisfy the Virasoro algebra of central charge $c$:

\begin{equation*}
\lbrack L(m),L(n)]=(m-n)L(m+n)+(m^{3}-m)\frac{c}{12}\delta _{m,-n}.
\end{equation*}
We define the homogeneous space of weight $k$ to be $V_{k}=\{v\in
V|L(0)v=kv\}$ where we write $\text{wt}(v)=k$ for $v$ in $V_{k}$ . Then as
an operator on $V$ we have

\begin{equation*}
v(n):V_{m}\rightarrow V_{m+k-n-1}.
\end{equation*}%
In particular, the \textit{zero mode} $o(v)=v(\text{wt}(v)-1)$ is a linear
operator on $V_{m}$. A state $v$ is said to be \textit{quasi-primary} if $%
L(1)v=0$ and \textit{primary} if additionally $L(2)v=0$.

The subalgebra $\{L(-1),L(0),L(1)\}$ generates a natural action on vertex
operators associated with $SL(2,\mathbb{C})$ M\"{o}bius transformations on $%
z $ (\cite{B1, DGM, FHL, Ka}). In particular, we note the inversion $%
z\mapsto 1/z$ for which 
\begin{equation}
Y(v,z)\mapsto Y^{\dagger }(v,z)=Y(e^{zL(1)}(-\frac{1}{z^{2}})^{L(0)}v,\frac{1%
}{z}).  \label{eq: adj op}
\end{equation}%
$Y^{\dagger }(v,z)$ is the \emph{adjoint} vertex operator \cite{FHL}. Under
the dilatation $z\mapsto az$ we have 
\begin{equation}
Y(v,z)\mapsto a^{L(0)}Y(v,z)a^{-L(0)}=Y(a^{L(0)}v,az).  \label{Y_D}
\end{equation}%
We also note (\cite{BPZ, Z2}) that under a general origin-preserving
conformal map $z\mapsto w=\phi (z)$, 
\begin{equation}
Y(v,z)\mapsto Y((\phi ^{\prime }(z))^{L(0)}v,w),  \label{Y_phi}
\end{equation}%
for any primary vector $v$.

\bigskip We consider some particular VOAs, namely Heisenberg free boson and
lattice VOAs. Consider an $l$-dimensional complex vector space (i.e.,
abelian Lie algebra) $\mathfrak{H}$ equipped with a non-degenerate,
symmetric, bilinear form $(\ ,)$ and a distinguished orthonormal basis $%
a_{1},a_{2},\ldots a_{l}$. The corresponding affine Lie algebra is the
Heisenberg Lie algebra $\mathfrak{\hat{H}}=\mathfrak{H}\otimes \mathbb{C}%
[t,t^{-1}]\oplus \mathbb{C}k$ with brackets $[k,\mathfrak{\hat{H}}]=0$ and

\begin{equation}
\lbrack a_{i}\otimes t^{m},a_{j}\otimes t^{n}]=m\delta _{i,j}\delta _{m,-n}k.
\label{Fockbracket}
\end{equation}%
Corresponding to an element $\lambda $ in the dual space $\mathfrak{H}^{\ast
}$ we consider the Fock space defined by the induced (Verma) module 
\begin{equation*}
M^{(\lambda )}=U(\mathfrak{\hat{H}})\otimes _{U(\mathfrak{H}\otimes \mathbb{C%
}[t]\oplus \mathbb{C}k)}\mathbb{C},
\end{equation*}%
where $\mathbb{C}$ is the $1$-dimensional space annihilated by $\mathfrak{H}%
\otimes t\mathbb{C}[t]$ and on which $k$ acts as the identity and $\mathfrak{%
H}\otimes t^{0}$ via the character $\lambda $; $U$ denotes the universal
enveloping algebra. There is a canonical identification of linear spaces

\begin{equation*}
M^{(\lambda )}=S(\mathfrak{H}\otimes t^{-1}\mathbb{C}[t^{-1}]),
\end{equation*}%
where $S$ denotes the (graded) symmetric algebra. The Heisenberg free boson
VOA $M^{l}$ corresponds to the case $\lambda =0$. The Fock states

\begin{equation}
v=a_{1}(-1)^{e_{1}}.a_{1}(-2)^{e_{2}}\ldots a_{1}(-n)^{e_{n}}\ldots
a_{l}(-1)^{f_{1}}.a_{l}(-2)^{f_{2}}\ldots a_{l}(-p)^{f_{p}}.\mathbf{1},
\label{Fockstate}
\end{equation}%
for non-negative integers $e_{i},\ldots ,f_{j}$ form a basis of $M^{l}$ with 
$a_i(n)\equiv a_i\otimes t^n$. The vacuum $\mathbf{1}$ is canonically
identified with the identity of $M_{0}=\mathbb{C}$, while the weight 1
subspace $M_{1}$ may be naturally identified with $\mathfrak{H}$. $M^{l}$ is
a simple VOA of central charge $l$.

\bigskip Next we consider the case of a lattice vertex operator algebra $%
V_{L}$ associated to a positive-definite even lattice $L$ (cf. \cite{B1, FLM}%
). Thus $L$ is a free abelian group of rank $l$ equipped with a positive
definite, integral bilinear form $(\ ,):L\otimes L\rightarrow \mathbb{Z}$
such that $(\alpha ,\alpha )$ is even for $\alpha \in L$. Let $\mathfrak{H}$
be the space $\mathbb{C}\otimes _{\mathbb{Z}}L$ equipped with the $\mathbb{C}
$-linear extension of $(\ ,)$ to $\mathfrak{H}\otimes \mathfrak{H}$ and let $%
M^{l}$ be the corresponding Heisenberg VOA. The Fock space of the lattice
theory may be described by the linear space 
\begin{equation}
V_{L}=M^{l}\otimes \mathbb{C}[L]=\sum_{\alpha \in L}M^{l}\otimes e^{\alpha },
\label{VLdefn}
\end{equation}%
where $\mathbb{C}[L]$ denotes the group algebra of $L$ with canonical basis $%
e^{\alpha }$, $\alpha \in L$. $M^{l}$ may be identified with the subspace $%
M^{l}\otimes e^{0}$ of $V_{L}$, in which case $M^{l}$ is a subVOA of $V_{L}$
and the rightmost equation of (\ref{VLdefn}) then displays the decomposition
of $V_{L}$ into irreducible $M^{l}$-modules. $V_{L}$ is a simple VOA of
central charge $l$. Each $\mathbf{1}\otimes e^{\alpha }\in V_{L}$ is a
primary state of weight $\frac{1}{2}(\alpha ,\alpha )$ with vertex operator
(loc. cit.) 
\begin{eqnarray}
Y(\mathbf{1}\otimes e^{\alpha },z) &=&Y_{-}(\mathbf{1}\otimes e^{\alpha
},z)Y_{+}(\mathbf{1}\otimes e^{\alpha },z)e^{\alpha }z^{\alpha },  \notag \\
Y_{\pm }(\mathbf{1}\otimes e^{\alpha },z) &=&\exp (\mp \sum_{n>0}\frac{%
\alpha (\pm n)}{n}z^{\mp n}).  \label{Yealpha}
\end{eqnarray}%
The operators $e^{\alpha }\in \mathbb{C}[L]$ obey 
\begin{equation}
e^{\alpha }e^{\beta }=\epsilon (\alpha ,\beta )e^{\alpha +\beta }
\label{eps_cocycle}
\end{equation}%
for $2$-cocycle $\epsilon (\alpha ,\beta )$ satisfying $\epsilon (\alpha
,\beta )\epsilon (\beta ,\alpha )=(-1)^{(\alpha ,\beta )}$.

\subsection{The Li-Zamolodchikov Metric}

\label{subsect_LiZ}

A bilinear form $\langle \ ,\rangle :V\times V{\longrightarrow }\mathbb{C}$
is called \emph{invariant} in case the following identity holds for all $%
a,b,c\in V$ (\cite{FHL}): 
\begin{equation}
\langle Y(a,z)b,c\rangle =\langle b,Y^{\dagger }(a,z)c\rangle,
\label{eq: inv bil form}
\end{equation}%
with $Y^{\dagger }(a,z)$ the adjoint operator (\ref{eq: adj op}).

\begin{remark}
\label{Rem_Zam}Note that 
\begin{eqnarray}
\langle a,b\rangle &=&\mathrm{Res}_{w=0}w^{-1}\mathrm{Res}%
_{z=0}z^{-1}\langle Y(a,w)\mathbf{1},Y(b,z)\mathbf{1}\rangle  \notag \\
&=&\mathrm{Res}_{w=0}w^{-1}\mathrm{Res}_{z=0}z^{-1}\langle \mathbf{1}%
,Y^{\dagger }(a,w)Y(b,z)\mathbf{1}\rangle  \notag \\
&=&"\langle \mathbf{1},Y(a,z=\infty )Y(b,z=0)\mathbf{1}\rangle " ,
\label{Zam}
\end{eqnarray}%
with $w=1/z$, following (\ref{eq: adj op}). Thus the invariant bilinear form
is equivalent to what is known as the (chiral) Zamolodchikov metric in
Conformal Field Theory (\cite{BPZ, P}).
\end{remark}

First note that any invariant bilinear form on $V$  is necessarily symmetric by a
theorem of \cite{FHL}. Generally a VOA may have \emph{no} non-zero invariant bilinear form, even if it is
well-behaved in other ways. Examples where $V$ is rational can be found in \cite{DM}. Results
of Li \cite{Li} guarantee that if $V_{0}$ is spanned by the vacuum vector $\mathbf{1}$ then the following hold: (a) $V$ has \emph{at most one} nonzero invariant bilinear form up to
scalars; (b) if $V$ \emph{has} a nonzero invariant bilinear form $\langle \ , \ \rangle$ then
the radical Rad$\langle \ , \ \rangle$ is the unique maximal ideal of $V$, and in particular $V$ is simple if,
and only if, $\langle \ , \ \rangle$ is
non-degenerate.  In this case, $V$ is \emph{self-dual} in the sense that $V$ is isomorphic to
the contragredient module $V^{\prime }$ as a $V$-module. Conversely, if $V$ is a self-dual VOA then it has a nondegenerate invariant bilinear form. All of the VOAs
that occur in this paper satisfy these conditions, i.e., they are simple and self-dual with $V_0 = \mathbb{C}\mathbf{1}$. Then if we normalize
so that $\langle \mathbf{1},\mathbf{1}\rangle =1$ then $\langle \ ,\rangle $
is unique and nondegenerate. We refer to this particular bilinear form as the \emph{%
Li-Zamolodchikov metric} on $V$, or LiZ-metric for short.

\begin{remark}
\label{Rem_tensorVOA}Uniqueness entails that the LiZ-metric on the tensor
product $V_{1}\otimes V_{2}$ of a pair of simple VOAs satisfying the
appropriate conditions is just the tensor product of the LiZ metrics on $%
V_{1}$ and $V_{2}$.
\end{remark}

If $a$ is a homogeneous, quasi-primary state, the component form of (\ref%
{eq: inv bil form}) reads 
\begin{equation}
\langle a(n)b,c\rangle =(-1)^{\text{wt}(a)}\langle b,a(2\text{wt}%
(a)-n-2)c\rangle .  \label{eq: qp bil form}
\end{equation}%
In particular, since the conformal vector $\omega $ is quasi-primary of
weight $2$ we may take $\omega $ in place of $a$ in (\ref{eq: qp bil form})
and obtain 
\begin{equation}
\langle L(n)b,c\rangle =\langle b,L(-n)c\rangle .
\label{eq: conf vec bil form}
\end{equation}%
The case $n=0$ of (\ref{eq: conf vec bil form}) shows that the homogeneous
spaces $V_{n},V_{m}$ are orthogonal if $n\not=m$. Taking $u=\mathbf{1}$ and
using $a=a(-1)\mathbf{1}$ in (\ref{eq: qp bil form}) yields 
\begin{equation}
\langle a,b\rangle =(-1)^{\text{wt}(a)}\langle \mathbf{1},a(2\text{wt}%
(a)-1)b\rangle ,  \label{eq: liz form}
\end{equation}%
for $a$ quasi-primary, and this affords a practical way to compute the
LiZ-metric.

\bigskip Consider the rank one Heisenberg (free boson) VOA $M = M^1$
generated by a weight one state $a$ with $(a,a)=1$. Then $\langle a,a\rangle
=-\langle \mathbf{1},a(1)a(-1)\mathbf{1}\rangle =-1$. Using (\ref%
{Fockbracket}), it is straightforward to verify that in general the Fock
basis consisting of vectors of the form 
\begin{equation}
v=a(-1)^{e_{1}}\ldots a(-p)^{e_{p}}.\mathbf{1,}  \label{eq: fock basis}
\end{equation}%
for non-negative integers $\{e_{i}\}$ is orthogonal with respect to the
LiZ-metric, and that 
\begin{equation}
\langle v,v\rangle =\prod_{1\leq i\leq p}(-i)^{e_{i}}e_{i}!.
\label{eq: inner prod}
\end{equation}%
This result generalizes in an obvious way for a rank $l$ free boson VOA $%
M^{l}$ with Fock basis (\ref{Fockstate}) following Remark \ref{Rem_tensorVOA}%
.

\section{Partition and $n$-Point Functions for Vertex Operator Algebras on a
Riemann Surface}

\label{sect_n point}

In this section we consider the partition and $n$-point functions for a VOA
on a Riemann surface of genus zero, one or two. Our definitions are based on
sewing schemes for the given Riemann surface in terms of one or more
surfaces of lower genus and are motivated by ideas in conformal field theory
especially \cite{FS, So1, P}. We assume throughout that $V$ has a
non-degenerate LiZ metric $\langle \ ,\rangle $. Then for any $V$ basis $%
\{u^{(a)}\}$, we may define the \emph{dual basis }$\{\bar{u}^{(a)}\}$ with
respect to the LiZ metric where 
\begin{equation}
\langle u^{(a)},\bar{u}^{(b)}\rangle =\delta _{ab}\text{.}  \label{LiZdual}
\end{equation}

\subsection{Genus Zero Case}

\label{subsect_genuszero}

We begin with the definition of the genus zero $n$-point function given by:%
\begin{equation}
Z_{V}^{(0)}(v_{1},z_{1};\ldots v_{n},z_{n})=\langle \mathbf{1}%
,Y(v_{1},z_{1})\ldots Y(v_{n},z_{n})\mathbf{1}\rangle ,  \label{Zzero}
\end{equation}%
for $v_{1},\ldots v_{n}\in V$. In particular, the genus zero partition (or $%
0 $-point) function is $Z_{V}^{(0)}=\langle \mathbf{1},\mathbf{1}\rangle =1$%
. The genus zero $n$-point function is a rational function of $z_{1},\ldots
z_{n}$, which we refer to as the insertion points, with possible poles at $%
z_{i}=0$ and $z_{i}=z_{j},i\neq j$. Thus we may consider $z_{1},\ldots
z_{n}\in \mathbb{C}\cup \{\infty \}$, the Riemann sphere, with $%
Z_{V}^{(0)}(v_{1},z_{1};\ldots ;v_{n},z_{n})$ evaluated for $\left\vert
z_{1}\right\vert >\left\vert z_{2}\right\vert >\ldots >\left\vert
z_{n}\right\vert $ (e.g. \cite{FHL, Z2, GG}). \ The $n$-point function has a
canonical geometric interpretation for primary vectors $v_{i}$ of $L(0)$
weight $\text{wt}(v_{i})$. Then $Z_{V}^{(0)}(v_{1},z_{1};\ldots
;v_{n},z_{n}) $ parameterizes a global meromorphic differential form on the
Riemann sphere, 
\begin{equation}
\mathcal{F}_{V}^{(0)}(v_{1},\ldots v_{n})=Z_{V}^{(0)}(v_{1},z_{1};\ldots
;v_{n},z_{n})\prod_{1\leq i\leq n}(dz_{i})^{\text{wt}(v_{i})}.
\label{Zzeroform}
\end{equation}%
It follows from (\ref{Y_phi}) that $\mathcal{F}_{V}^{(0)}$ is conformally
invariant.\ This interpretation is the starting point of various
algebraic-geometric approaches to $n$-point functions (\emph{apart }from the
partition or $0$-point function) at higher genera (e.g. \cite{TUY, Z2}).

\medskip It is instructive to consider $\mathcal{F}_{V}^{(0)}$ in the
context of a trivial sewing of two Riemann spheres parameterized by $z_{1}$
and $z_{2}$ to form another Riemann sphere as follows. For $r_{a}>0,a=1,2$,
and a complex parameter $\epsilon $ satisfying $|\epsilon |\leq r_{1}r_{2}$,
excise the open disks $\left\vert z_{a}\right\vert <|\epsilon |r_{\bar{a}%
}^{-1}$ (recall convention (\ref{abar})) and identify the annular regions $%
r_{a}\geq \left\vert z_{a}\right\vert \geq |\epsilon |r_{\bar{a}}^{-1}$ via
the sewing relation 
\begin{equation}
z_{1}z_{2}=\epsilon.  \label{spheresew}
\end{equation}%
Consider $Z_{V}^{(0)}(v_{1},x_{1};\ldots v_{n},x_{n})$ for quasi-primary $%
v_{i}$ with $r_{1}\geq |x_{i}|\geq |\epsilon |r_{2}^{-1}$ and let $%
y_{i}=\epsilon /x_{i}$. Then for $0\leq k\leq n-1$ we find from (\ref%
{LiZdual}) that 
\begin{eqnarray}
&&\hspace{1cm}Y(v_{k+1},x_{k+1})\ldots Y(v_{n},x_{n})\mathbf{1}=  \notag \\
&&\sum_{r\geq 0}\sum_{u\in V_{r}}\langle \bar{u},Y(v_{k+1},x_{k+1})\ldots
Y(v_{n},x_{n})\mathbf{1}\rangle u,  \label{Yxkplus_u_expansion}
\end{eqnarray}%
where the inner sum is taken over any basis for $V_{r}$. Thus 
\begin{eqnarray*}
&&\hspace{3.5cm}Z_{V}^{(0)}(v_{1},x_{1};\ldots v_{n},x_{n})= \\
&&\sum_{r\geq 0}\sum_{u\in V_{r}}\langle \mathbf{1},Y(v_{1},x_{1})\ldots
Y(v_{k},x_{k})u\rangle \langle \bar{u},Y(v_{k+1},x_{k+1})\ldots
Y(v_{n},x_{n})\mathbf{1}\rangle .
\end{eqnarray*}%
But 
\begin{equation*}
\langle \mathbf{1},Y(v_{1},x_{1})\ldots Y(v_{k},x_{k})u\rangle =\mathrm{Res}%
_{z_{1}=0}z_{1}^{-1}Z_{V}^{(0)}(v_{1},x_{1};\ldots v_{k},x_{k};u,z_{1}),
\end{equation*}%
and 
\begin{eqnarray*}
&&\langle \bar{u},Y(v_{k+1},x_{k+1})\ldots Y(v_{n},x_{n})\mathbf{1}\rangle \\
&=&\langle \mathbf{1},Y^{\dagger }(v_{n},x_{n})\ldots Y^{\dagger
}(v_{k+1},x_{k+1})\bar{u}\rangle \\
&=&\langle \mathbf{1},\epsilon ^{L(0)}Y^{\dagger }(v_{n},x_{n})\epsilon
^{-L(0)}\ldots \epsilon ^{L(0)}Y^{\dagger }(v_{k+1},x_{k+1})\epsilon
^{-L(0)}\epsilon ^{L(0)}\bar{u}\rangle \\
&=&\epsilon ^{r}\mathrm{Res}_{z_{2}=0}z_{2}^{-1}Z_{V}^{(0)}(v_{n},y_{n};%
\ldots v_{k+1},y_{k+1};\bar{u},z_{2})\prod_{k+1\leq i\leq n}(-\frac{\epsilon 
}{x_{i}^{2}})^{\text{wt}(v_{i})}.
\end{eqnarray*}%
The last equation holds since for quasiprimary states $v_{i}$, the M\"{o}%
bius transformation $x\mapsto y=\epsilon /x$ induces 
\begin{equation}
Y(v_{i},x_{i})\mapsto \epsilon ^{L(0)}Y^{\dagger }(v_{i},x_{i})\epsilon
^{-L(0)}=(-\frac{\epsilon }{x_{i}^{2}})^{\text{wt}(v_{i})}Y(v_{i},y_{i}).
\label{YMobius}
\end{equation}%
Thus we find:

\begin{proposition}
\label{Prop_Z0sew} For homogeneous quasiprimary states $v_{i}$ with the
sewing scheme (\ref{spheresew}), we have 
\begin{gather*}
\mathcal{F}_{V}^{(0)}(v_{1},\ldots ,v_{n})= \\
\sum_{r\geq 0}\epsilon ^{r}\sum_{u\in V_{r}}\mathrm{Res}%
_{z_{1}=0}z_{1}^{-1}Z_{V}^{(0)}(v_{1},x_{1};\ldots v_{k},x_{k};u,z_{1}) \\
\mathrm{Res}_{z_{2}=0}z_{2}^{-1}Z_{V}^{(0)}(v_{n},y_{n};\ldots
v_{k+1},y_{k+1};\bar{u},z_{2})\prod_{1\leq i\leq k}(dx_{i})^{\text{wt}%
(v_{i})}\prod_{k+1\leq i\leq n}(dy_{i})^{\text{wt}(v_{i})},
\end{gather*}
for any $k$, $0\leq k\leq n-1$ i.e. the RHS is independent of the choice of
Riemann sphere on which the insertion point of each state $v_i$ lies. $%
\square$
\end{proposition}

\subsection{Genus One Case}

We now consider genus one $n$-point functions defined in terms of a
self-sewing of a Riemann sphere where punctures are located at the origin
and the point at infinity \cite{MT2}. Choose local coordinates $z_{1}=z$ in
the neighborhood of the origin and $z_{2}=1/z^{\prime }$ for $z^{\prime }$
in the neighborhood of the point at infinity. For $a=1,2$ and $r_{a}>0$,
identify the annular regions $|q|r_{\bar{a}}^{-1}\leq \left\vert
z_{a}\right\vert \leq r_{a}$ for complex $q$ satisfying $|q|\leq r_{1}r_{2}$
via the sewing relation $z_{1}z_{2}=q$ i.e. $z=qz^{\prime }$. Then it is
straightforward to show that the annuli do not intersect for $|q|<1$, and
that $q=\exp (2\pi i\tau )$ where $\tau $ is the torus modular parameter
(e.g. \cite{MT2}, Proposition 8).

\medskip We define the genus one partition function by 
\begin{eqnarray}
&&Z_{V}^{(1)}(q)=Z_{V}^{(1)}(\tau )=  \notag \\
&&q^{-c/24} \sum_{n\geq 0}q^{n}\sum_{u\in V_{n}}\mathrm{Res}%
_{z_{2}=0}z_{2}^{-1}\mathrm{Res}_{z_{1}=0}z_{1}^{-1}\langle \mathbf{1}%
,Y^{\dag }(u,z_{2})Y(\bar{u},z_{1})\mathbf{1}\rangle ,  \label{Z1def}
\end{eqnarray}%
where the inner sum is taken over any basis for $V_{n}$. The external factor
of $q^{-c/24}$ is introduced in the usual way to enhance the modular
properties of $Z_{V}^{(1)}(q)$ \cite{Z1}. From (\ref{Zam}) and (\ref{LiZdual}%
) it follows that 
\begin{equation}
Z_{V}^{(1)}(\tau )=\sum_{n\geq 0}\dim V_{n}q^{n-c/24}=\mathrm{Tr}%
_{V}(q^{L(0)-c/24}),  \label{Z_tau}
\end{equation}%
the standard graded trace definition. The genus one $n$-point function is
similarly given by%
\begin{eqnarray*}
&&\sum_{r\geq 0}q^{r-c/24}\sum_{u\in V_{r}}\mathrm{Res}_{z_{2}=0}z_{2}^{-1}%
\mathrm{Res}_{z_{1}=0}z_{1}^{-1}\langle \mathbf{1},Y^{\dag
}(u,z_{2})Y(v_{1},x_{1})\ldots Y(v_{n},x_{n})Y(\bar{u},z_{1})\mathbf{1}%
\rangle \\
&&=\mathrm{Tr}_{V}(Y(v_{1},x_{1})\ldots Y(v_{n},x_{n})q^{L(0)-c/24}).
\end{eqnarray*}%
\qquad It is natural to consider the conformal map $x=q_{z}\equiv \exp (z)$
in order to describe the elliptic properties of the $n$-point function \cite%
{Z1}. Since from (\ref{Y_phi}), for a primary state $v$, $Y(v,w)\rightarrow
Y(q_{z}^{L(0)}v,q_{z})$ under this conformal map, we are led to the
following definition of the genus one $n$-point function (op. cite.): 
\begin{eqnarray}
&&Z_{V}^{(1)}(v_{1},z_{1};\ldots v_{n},z_{n};\tau )=  \notag \\
&&\mathrm{Tr}_{V}(Y(q_{z_{1}}^{L(0)}v_{1},q_{z_{1}})\ldots
Y(q_{z_{n}}^{L(0)}v_{n},q_{z_{n}})q^{L(0)-c/24}).  \label{Z1npt}
\end{eqnarray}%
For homogeneous primary states $v_{i}$ of weight $\text{wt}(v_{i})$, $%
Z_{V}^{(1)}$ parameterizes a global meromorphic differential form on the
torus 
\begin{equation}
\mathcal{F}_{V}^{(1)}(v_{1},\ldots v_{n};\tau
)=Z_{V}^{(1)}(v_{1},z_{1};\ldots v_{n},z_{n};\tau )\prod_{1\leq i\leq
n}(dz_{i})^{\text{wt}(v_{i})}.  \label{F1npt}
\end{equation}

\bigskip Zhu introduced (\cite{Z1}) a second VOA $(V,Y[,],\mathbf{1},\tilde{%
\omega})$ which is isomorphic to $(V,Y(,),\mathbf{1},\omega ).$ It has
vertex operators 
\begin{equation}
Y[v,z]=\sum_{n\in \mathbb{Z}}v[n]z^{-n-1}=Y(q_{z}^{L(0)}v,q_{z}-1),
\label{Ysquare}
\end{equation}%
and conformal vector $\tilde{\omega}=$ $\omega -\frac{c}{24}\mathbf{1}$. Let 
\begin{equation}
Y[\tilde{\omega},z]=\sum_{n\in \mathbb{Z}}L[n]z^{-n-2} ,  \label{Ywtilde}
\end{equation}%
and write $\text{wt}[v]=k$ if $L[0]v=kv$, $V_{[k]}=\{v\in V|\text{wt}[v]=k\}$%
. Only primary vectors are homogeneous with respect to both $L(0)$ and $L[0]$%
, in which case $\text{wt}(v)=\text{wt}[v]$. Similarly, we define the square
bracket LiZ metric $\langle \ ,\rangle _{\mathrm{sq}}$ which is invariant
with respect to the square bracket adjoint.

\medskip We denote $1$-point functions by 
\begin{equation}
Z_{V}^{(1)}(v,\tau )=Z_{V}^{(1)}(v,z;\tau )=\mathrm{Tr}%
_{V}(o(v)q^{L(0)-c/24}).  \label{Z1_1_pt}
\end{equation}%
($Z_{V}^{(1)}(v,\tau )$ is necessarily $z$ independent). Any $n$-point
function can be expressed in terms of $1$-point functions (\cite{MT1}, Lemma
3.1) as follows: 
\begin{eqnarray}
&&Z_{V}^{(1)}(v_{1},z_{1};\ldots v_{n},z_{n};\tau )  \notag \\
&=&Z_{V}^{(1)}(Y[v_{1},z_{1}]\ldots Y[v_{n-1},z_{n-1}]Y[v_{n},z_{n}]\mathbf{1%
},\tau )  \label{Z1Ysq1} \\
&=&Z_{V}^{(1)}(Y[v_{1},z_{1n}]\ldots Y[v_{n-1},z_{n-1n}]v_{n},\tau ),
\label{Z1Ysq2}
\end{eqnarray}%
where $z_{in}=z_{i}-z_{n}$.

We may consider a trivial sewing of a torus with local coordinate $z_{1}$ to
a Riemann sphere with local coordinate $z_{2}$ by identifying the annuli $%
r_{a}\geq |z_{a}|\geq |\epsilon |r_{\bar{a}}^{-1}$ via the sewing relation $%
z_{1}z_{2}=\epsilon$. Consider $Z_{V}^{(1)}(v_{1},x_{1};\ldots v_{n},x_{n})$
for quasi-primary $v_{i}$ of $L[0]$ weight $\text{wt}[v_{i}]$, with $%
r_{1}\geq |x_{i}|\geq |\epsilon |r_{2}^{-1}$, and let $y_{i}=\epsilon /x_{i}$%
. Using (\ref{Z1Ysq1}), and employing the square bracket version of (\ref%
{Yxkplus_u_expansion}) with square bracket LiZ metric $\langle \ ,\rangle _{%
\mathrm{sq}}$, we have 
\begin{eqnarray*}
&&Z_{V}^{(1)}(v_{1},x_{1};\ldots v_{n},x_{n};\tau )= \\
&&\sum_{r\geq 0}\sum_{u\in V_{[r]}}Z_{V}^{(1)}(Y[v_{1},x_{1}]\ldots
Y[v_{k},x_{k}]u;\tau )\langle \bar{u},Y[v_{k+1},x_{k+1}]\ldots Y[v_{n},x_{n}]%
\mathbf{1}\rangle _{\mathrm{sq}},
\end{eqnarray*}%
where the inner sum is taken over any basis $\{u\}$ of $V_{[r]}$, and $\{%
\bar{u}\}$ is the dual basis with respect to $\langle \ ,\rangle _{\mathrm{sq%
}}$. Now 
\begin{equation*}
Z_{V}^{(1)}(Y[v_{1},x_{1}]\ldots Y[v_{k},x_{k}]u;\tau )=\mathrm{Res}%
_{z_{1}=0}z_{1}^{-1}Z_{V}^{(1)}(v_{1},x_{1};\ldots v_{k},x_{k};u,z_{1};\tau).
\end{equation*}%
Using the isomorphism between the round and square bracket formalisms, we
find as before that 
\begin{eqnarray*}
&&\langle \bar{u},Y[v_{k+1},x_{k+1}]\ldots Y[v_{n},x_{n}]\mathbf{1}\rangle _{%
\mathrm{sq}} \\
&=&\epsilon ^{r}\mathrm{Res}_{z_{2}=0}z_{2}^{-1}Z_{V}^{(0)}(v_{n},y_{n};%
\ldots v_{k+1},y_{k+1};\bar{u},z_{2}) \prod_{k+1\leq i\leq n}(-\frac{%
\epsilon }{x_{i}^{2}})^{\text{wt}[v_{i}]}.
\end{eqnarray*}%
We thus obtain a natural analogue of Proposition \ref{Prop_Z0sew}:

\begin{proposition}
\label{Prop_Z1sew} For square bracket homogeneous quasiprimary states $v_{i}$
with the above sewing scheme, then we have 
\begin{gather*}
\mathcal{F}_{V}^{(1)}(v_{1},\ldots ,v_{n};\tau )= \\
\sum_{r\geq 0}\epsilon ^{r}\sum_{u\in V_{[r]}}\mathrm{Res}%
_{z_{1}=0}z_{1}^{-1}Z_{V}^{(1)}(v_{1},x_{1};\ldots v_{k},x_{k};u,z_{1};\tau).
\\
\mathrm{Res}_{z_{2}=0}z_{2}^{-1}Z_{V}^{(0)}(v_{n},y_{n};\ldots
v_{k+1},y_{k+1};\bar{u},z_{2}) \prod_{1\leq i\leq k}(dx_{i})^{\text{wt}%
[v_{i}]} \prod_{k+1\leq i\leq n}(dy_{i})^{\text{wt}[v_{i}]},\quad
\end{gather*}%
and is independent of $k=0,1,\ldots n-1$ where the inner sum is taken over
any basis $\{u\}$ for $V_{[r]}$, $\{\bar{u}\}$ is the dual basis with
respect to $\langle \ ,\rangle _{\mathrm{sq}}$. \hfill $\square $
\end{proposition}

We note that all the above definitions can be naturally extended for any $V$%
-module $N$ with vertex operators $Y_N(v,x)$ where the trace in (\ref{Z1Ysq2}%
) is taken over $N$ and $o(v)$ is replaced by $o_N(v)$ the Virasoro level
preserving part of $Y_N(v,x)$.

\subsection{Genus Two Case}

Motivated by Proposition \ref{Prop_Z1sew}, we now discuss the formal
definition of the genus two $n$-point function associated with the genus two 
$\epsilon $-sewing scheme reviewed in Section \ref{subsect_epsilon}. Recall
that we sew together a pair of tori $\mathcal{S}_{1},\mathcal{S}_{2}$ with
modular parameters $\tau _{1},\tau_{2}$ respectively via the sewing relation
(\ref{pinch}). We define the genus two $n$-point function for $v_1,\ldots
v_k $ inserted at $x_{1},\ldots ,x_{k}\in \hat{\mathcal{S}_{1}}$ and $%
v_{k+1},\ldots v_n$ inserted at $y_{k+1},\ldots ,y_{n}\in \hat{\mathcal{S}%
_{2}}$ for $k=0,1,\ldots n-1$ by 
\begin{eqnarray}
&&Z_{V}^{(2)}(v_{1},x_{1};\ldots v_{k},x_{k}\vert v_{k+1},y_{k+1};\ldots
v_{n},y_{n};\tau _{1},\tau _{2},\epsilon )  \notag \\
&&=\sum_{r\geq 0}\epsilon ^{r}\sum_{u\in V_{[r]}}\mathrm{Res}%
_{z_{1}=0}z_{1}^{-1} Z_{V}^{(1)}(v_{1},x_{1};\ldots
v_{k},x_{k};u,z_{1};\tau_{1})\cdot  \notag \\
&&\ \ \ \mathrm{Res}_{z_{2}=0}z_{2}^{-1} Z_{V}^{(1)}(v_{n},y_{n};\ldots
v_{k+1},y_{k+1};\bar{u},z_{2}; \tau _{2}),  \notag \\
&&=\sum_{r\geq 0}\epsilon ^{r}\sum_{u\in V_{[r]}}
Z_{V}^{(1)}(Y[v_{1},x_{1}]\ldots Y[v_{k},x_{k}]u,z_{1};\tau_{1})\cdot  \notag
\\
&&\ \ \ Z_{V}^{(1)}(Y[v_{n},y_{n}]\ldots Y[v_{k+1},y_{k+1}]\bar{u},z_{2};
\tau _{2}),  \label{Z2n_pt_eps}
\end{eqnarray}%
where the inner sum is taken over any basis $V_{[r]}$ and $\bar{u}$ is the
dual of $u$ with respect to $\langle \ ,\rangle _{\mathrm{sq}}$. The last
expression in (\ref{Z2n_pt_eps}) follows from (\ref{Z1Ysq2}).

\begin{remark}
\label{Rem_tensorVOAZ2} Following Remark \ref{Rem_tensorVOA} it is clear
that the genus two $n$-point function on the tensor product $V_{1}\otimes
V_{2}$ of a pair of simple VOAs is just the product of $n$-point functions
on $V_{1}$ and $V_{2}$.
\end{remark}

In this paper we mainly concentrate on the genus two partition function
(i.e. the $0$-point function) given by 
\begin{equation}
Z_{V}^{(2)}(\tau _{1},\tau _{2},\epsilon )=\sum_{n\geq 0}\epsilon
^{n}\sum_{u\in V_{[n]}}Z_{V}^{(1)}(u,\tau _{1})Z_{V}^{(1)}(\bar{u},\tau
_{2}).  \label{Z2_def_eps}
\end{equation}
Some examples of $n$-point functions will also be computed. A general
discussion of all genus two $n$-point functions for the Heisenberg VOA and
its modules will appear elsewhere \cite{MT5}.

Clearly the definition of the $n$-point function (\ref{Z2n_pt_eps}) depends
on the choice of punctured torus on which the insertion points lie. However,
by defining an associated formal differential form, we find the following
genus two analogue of Propositions \ref{Prop_Z0sew} and \ref{Prop_Z1sew}:

\begin{proposition}
\label{Prop_Z2sew} For $x_{i}\in \hat{\mathcal{S}_{1}}$ and $y_i\in \hat{%
\mathcal{S}_{2}}$ with $x_i y_i=\epsilon$ and square bracket homogeneous
quasiprimary states $v_{i}$, the formal differential form 
\begin{eqnarray}
&&\mathcal{F}_{V}^{(2)}(v_{1},\ldots ,v_{n};\tau _{1},\tau _{2},\epsilon
)\equiv  \notag \\
&&Z_{V}^{(2)}(v_{1},x_{1};\ldots v_{k},x_{k}\vert v_{k+1},y_{k+1};\ldots
v_{n},y_{n};\tau _{1},\tau _{2},\epsilon ) \cdot  \notag \\
&& \prod_{1\leq i\leq k}(dx_{i})^{\text{wt}[v_{i}]} \prod_{k+1\leq i\leq
n}(dy_{i})^{\text{wt}[v_{i}]},  \label{F2}
\end{eqnarray}
is independent of $k=0,1,\ldots n-1$.
\end{proposition}

\noindent \textbf{Proof.} Consider the left torus contribution in the
summand of (\ref{Z2n_pt_eps}) and expand $Y[v_{k},x_{k}]u$ in a square
bracket homogeneous basis: 
\begin{eqnarray*}
&& Z_{V}^{(1)}(Y[v_{1},x_{1}]\ldots Y[v_{k},x_{k}]u;\tau_{1})= \\
&&\sum_{s\ge 0} \sum_{w \in V_{[s]}}Z_{V}^{(1)}(Y[v_{1},x_{1}]\ldots
Y[v_{k-1},x_{k-1}]w;\tau_{1}) \langle \bar{w},Y[v_k,x_k] u\rangle_{\mathrm{sq%
}}.
\end{eqnarray*}
But for quasi-primary $v_k$ and using (\ref{YMobius}) we find 
\begin{eqnarray*}
\epsilon^r\langle \bar{w},Y[v_k,x_k] u\rangle_{\mathrm{sq}}= \langle
\epsilon^{L[0]}Y^\dagger[v_k,x_k]\bar{w},u\rangle_{\mathrm{sq}} = \epsilon^s
(-\frac{\epsilon }{x_{k}^{2}})^{\text{wt}[v_{k}]} \langle Y[v_k,y_k]\bar{w}%
,u\rangle_{\mathrm{sq}},
\end{eqnarray*}
where $x_ky_k=\epsilon$. Noting that 
\begin{eqnarray*}
&&\sum_{r\ge 0}\sum_{u \in V_{[r]}}Z_{V}^{(1)}(Y[v_{n},y_{n}]\ldots
Y[v_{k+1},y_{k+1}]\bar{u};\tau_{2}) \langle u,Y[v_k,y_k]\bar{w}\rangle_{%
\mathrm{sq}} \\
&&= Z_{V}^{(1)}(Y[v_{n},y_{n}]\ldots Y[v_{k+1},y_{k+1}]Y[v_k,y_k]\bar{w}%
;\tau_{2}),
\end{eqnarray*}
we therefore find that 
\begin{eqnarray*}
&&Z_{V}^{(2)}(v_{1},x_{1};\ldots v_{k},x_{k}\vert v_{k+1},y_{k+1};\ldots
v_{n},y_{n} ;\tau _{1},\tau _{2},\epsilon ) \\
=&&(-\frac{\epsilon }{x_{k}^{2}})^{\text{wt}[v_{k}]}
Z_{V}^{(2)}(v_{1},x_{1};\ldots v_{k-1},x_{k-1}\vert v_{k},y_{k};\ldots
v_{n},y_{n};\tau _{1},\tau _{2},\epsilon ).
\end{eqnarray*}
Hence 
\begin{eqnarray*}
&&Z_{V}^{(2)}(v_{1},x_{1};\ldots v_{k},x_{k}\vert v_{k+1},y_{k+1};\ldots
v_{n},y_{n};\tau _{1},\tau _{2},\epsilon ) \cdot  \notag \\
&& \prod_{1\leq i\leq k}(dx_{i})^{\text{wt}[v_{i}]} \prod_{k+1\leq i\leq
n}(dy_{i})^{\text{wt}[v_{i}]}= \\
&&Z_{V}^{(2)}(v_{1},x_{1};\ldots v_{k-1},x_{k-1}\vert v_{k},y_{k};\ldots
v_{n},y_{n};\tau _{1},\tau _{2},\epsilon ) \cdot  \notag \\
&& \prod_{1\leq i\leq k-1}(dx_{i})^{\text{wt}[v_{i}]} \prod_{k\leq i\leq
n}(dy_{i})^{\text{wt}[v_{i}]}.
\end{eqnarray*}
The result follows by repeated application of this identity. \hfill $\square$

\begin{remark}
\label{rem:F2global} If $Z_V^{(2)}(\tau _{1},\tau _{2},\epsilon )$ is
convergent on $\mathcal{D}^\epsilon$, we conjecture that for primary states $%
v_1,\ldots v_n$ then $\mathcal{F}_{V}^{(2)}(v_{1},\ldots ,v_{n};\tau
_{1},\tau _{2},\epsilon )$ is a genus two global meromorphic form with
possible poles only at coincident insertion points.
\end{remark}

Finally, note that all the above definitions can be naturally extended for
any pair of $V$-modules $N_1,N_2$ where the left (right) 1-point function in
(\ref{Z2n_pt_eps}) is considered for $N_1$ (respectively $N_2$).

\section{The Heisenberg VOA}

\label{sect_Boson} In this section we compute closed formulas for the genus
two partition function for the rank one Heisenberg VOA $M$ and compute the $%
n $-point function for $n$ Heisenberg vectors and the Virasoro vector
1-point function. We also discuss the modular properties of the partition
function in some detail.

\subsection{The Genus Two Partition Function $Z_{M}^{(2)}(\protect\tau _{1},%
\protect\tau _{2},\protect\epsilon )$}

\label{subsect_BosonZ2}

We wish to establish a closed formula for the genus two partition function $%
Z_{M}^{(2)}(\tau _{1},\tau _{2},\epsilon )$ of (\ref{Z2_def_eps}) in terms
of the infinite matrices $A_1,A_2$ introduced in (\ref{Akldef}) of Section
2. Recalling the definition (\ref{logdet}) we have:

\begin{theorem}
\label{Theorem_Z2_boson} Let $M$ be the vertex operator algebra of one free
boson. Then 
\begin{equation}
Z_{M}^{(2)}(\tau _{1},\tau _{2},\epsilon )=Z_{M}^{(1)}(\tau
_{1})Z_{M}^{(1)}(\tau _{2})(\det (I-A_{1}A_{2}))^{-1/2},  \label{Z2_1bos}
\end{equation}%
where $Z_{M}^{(1)}(\tau )=1/\eta (\tau )$.
\end{theorem}

\begin{remark}
\label{Rem_Z2mult} From Remark \ref{Rem_tensorVOAZ2} it follows that the
genus two partition function for $l$ free bosons $M^{l}$ is just the $l^{th}$
power of (\ref{Z2_1bos}).
\end{remark}

\noindent \textbf{Proof of Theorem.} The genus two partition function $%
Z_{M}^{(2)}(\tau _{1},\tau _{2},\epsilon )$ of (\ref{Z2_def_eps}) is $V$
basis independent. We choose the standard Fock vectors (in the square
bracket formulation) 
\begin{equation}
v=a[-1]^{e_{1}}\ldots a[-p]^{e_{p}}\mathbf{1}.  \label{eq: sq fock vec}
\end{equation}
Of course, these Fock vectors correspond in a natural 1-1 manner with
unrestricted partitions, the state $v$ (\ref{eq: sq fock vec}) corresponding
to a partition $\lambda =\{1^{e_{1}}\ldots p^{e_{p}}\}$ with $\vert
\lambda\vert=\sum_i e_{i}$ elements of $n=\sum_{1\leq i\leq p}ie_{i}$. We
sometimes write $v=v(\lambda )$ to indicate this correspondence.
Furthermore, following (\ref{eq: inner prod}) $v(\lambda )=(-1)^{\vert
\lambda\vert}\big(\prod_{1\leq i\leq p}i^{e_{i}}e_{i}!\big)\bar v(\lambda )$%
. Thus with this diagonal basis we have 
\begin{equation}
Z_{M}^{(2)}(\tau _{1},\tau _{2},\epsilon )= \sum_{\lambda =\{i^{e_{i}}\}}%
\frac{(-1)^{\vert\lambda\vert}}{\prod_{i}i^{e_{i}}e_{i}!}\epsilon ^{\sum
ie_{i}}Z_{M}^{(1)}(v(\lambda ),\tau_1 )Z_{M}^{(1)}(v(\lambda ),\tau_2 ).
\label{eq: part func diag}
\end{equation}

As discussed at length in \cite{MT1}, the partition $\lambda$ may be thought
of as a labeled set $\Phi =\Phi _{\lambda }$ with $e_{i}$ elements labeled $%
i $. One of the main results of \cite{MT1} (loc.cit. Corollary 1 and
eqn.(53)) is that for even $\vert\lambda\vert$ 
\begin{equation}
Z_{M}^{(1)}(v(\lambda ),\tau )=Z_{M}^{(1)}(\tau )\sum_{\phi \in F(\Phi
_{\lambda })}\Gamma (\phi ),  \label{eq: genus1 ptn func}
\end{equation}%
with 
\begin{equation}
\Gamma (\phi ,\tau )=\Gamma (\phi )=\prod_{(r,s)}C(r,s,\tau ),
\label{eq: Gamma}
\end{equation}%
for $C$ of (\ref{Ckldef}), where $\phi $ ranges over the elements of $F(\Phi
_{\lambda })$ (the fixed-point-free involutions in $\Sigma (\Phi _{\lambda
}) $) and $(r,s)$ ranges over the orbits of $\phi $ on $\Phi _{\lambda }$.
If $\vert\lambda\vert$ is odd then $Z_{M}^{(1)}(v(\lambda ),\tau )=0$. \medskip
With this notation, (\ref{eq: part func diag}) reads 
\begin{equation}
Z_{M}^{(2)}(\tau _{1},\tau _{2},\epsilon )=Z_{M}^{(1)}(\tau
_{1})Z_{M}^{(1)}(\tau _{2})\sum_{\lambda =\{i^{e_{i}}\}}\frac{E(\lambda )}{%
\prod_{i}i^{e_{i}}e_{i}!}\epsilon ^{\sum ie_{i}},
\label{eq: part func lambda}
\end{equation}%
where $\lambda $ ranges over all even $\vert\lambda\vert$ unrestricted
partitions and where we have set 
\begin{eqnarray}
E(\lambda ) &=&\sum_{\phi ,\psi \in F(\Phi _{\lambda })}\Gamma _{1}(\phi
)\Gamma _{2}(\psi ),  \label{eq: expr E} \\
\Gamma _{i}(\phi ) &=&\Gamma (\phi ,\tau _{i}).  \label{eq: expr Gamma}
\end{eqnarray}

We now analyze the nature of the expression $E(\lambda )$ more closely. This
will lead us to the connection between $Z^{(2)}(\tau _{1},\tau _{2},\epsilon
)$ and the chequered cycles discussed in Section 3.1. The idea is to use the
technique employed in the proof of Proposition 3.10 of \cite{MT1}. If we fix
for a moment a partition $\lambda $ then a pair of fixed-point-free
involutions $\phi ,\psi $ correspond (loc.cit.) to a pair of complete
matchings $\mu _{\phi },\mu _{\psi }$ on the labeled set $\Phi _{\lambda }$
which we may represent pictorially as

\begin{equation*}
\begin{array}{ccccc}
\overset{r_{1}}{\bullet } & \overset{1}{\longrightarrow } & \overset{s_{1}}{%
\bullet } & \overset{2}{\longrightarrow } & \overset{t_{1}}{\bullet } \\ 
\overset{r_{2}}{\bullet } & \overset{1}{\longrightarrow } & \overset{s_{2}}{%
\bullet } & \overset{2}{\longrightarrow } & \overset{t_{2}}{\bullet } \\ 
\vdots &  & \vdots & \vdots & \vdots \\ 
\overset{r_{b}}{\bullet } & \overset{1}{\longrightarrow } & \overset{s_{b}}{%
\bullet } & \overset{2}{\longrightarrow } & \overset{t_{b}}{\bullet }%
\end{array}%
\end{equation*}

\begin{center}
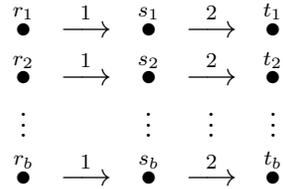
{\small Fig.~4 Two complete matchings }
\end{center}

\noindent Here, $\mu _{\phi }$ is the matching with edges labeled $1$, $\mu
_{\psi }$ the matching with edges labeled $2$, and where we have denoted the
(labeled) elements of $\Phi _{\lambda }$ by $\{r_{1},s_{1},\ldots
,r_{b},s_{b}\}=\{s_{1},t_{1},\ldots ,s_{b},t_{b}\}$. From this data we may
create a chequered cycle in a natural way: starting with some node of $\Phi
_{\lambda }$, apply the involutions $\phi ,\psi $ successively and
repeatedly until the initial node is reached, using the complete matchings
to generate a chequered cycle. The resulting chequered cycle corresponds to
an orbit of $\langle \psi \phi \rangle $ considered as a cyclic subgroup of $%
\Sigma (\Phi _{\lambda })$. Repeat this process for each such orbit to
obtain a \emph{chequered diagram} $D$ consisting of the union of the
chequered cycles corresponding to all of the orbits of $\langle \psi \phi
\rangle $ on $\Phi _{\lambda }$. To illustrate, for the partition $\lambda
=\{1^{2}.2.3^{2}.5\}$ with matchings $\mu _{\phi }=(13)(15)(23)$ and $\mu
_{\psi }=(11)(35)(23)$, the corresponding chequered diagram is

\begin{center}
\begin{picture}(300,60)

\put(100,50){\line(1,0){20}}
\put(95,53){\makebox(0,0){$1$}}
\put(100,50){\circle*{4}}
\put(110,54){\makebox(0,0){\scriptsize $1$}}

\put(120,50){\line(0,-1){20}}
\put(125,53){\makebox(0,0){$3$}}
\put(120,50){\circle*{4}}
\put(123,40){\makebox(0,0){\scriptsize $2$}}

\put(100,30){\line(1,0){20}}
\put(125,27){\makebox(0,0){$5$}}
\put(120,30){\circle*{4}}
\put(110,26){\makebox(0,0){\scriptsize $1$}}

\put(100,50){\line(0,-1){20}}
\put(95,27){\makebox(0,0){$1$}}
\put(100,30){\circle*{4}}
\put(97,40){\makebox(0,0){\scriptsize $2$}}

\put(160,40){\qbezier(0,0)(10,10)(20,0)}
\put(160,40){\circle*{4}}
\put(153,40){\makebox(0,0){$2$}}
\put(170,50){\makebox(0,0){\scriptsize $1$}}

\put(180,40){\circle*{4}}
\put(185,40){\makebox(0,0){$3$}}
\put(160,40){\qbezier(0,0)(10,-10)(20,0)}
\put(170,31){\makebox(0,0){\scriptsize $2$}}

\end{picture}

{\small Fig.~5 Chequered diagram}
\end{center}

Two chequered diagrams are isomorphic if there is a bijection on
the nodes which preserves edges and labels of nodes and edges. If $\lambda=\{1^{e_{1}}\ldots p^{e_{p}}\}$ then $\Sigma (\Phi _{\lambda })$ acts on the
chequered diagrams which have $\Phi _{\lambda }$ as underlying set of
labeled nodes. The \emph{Automorphism subgroup} $\mathrm{Aut}(D)$, consisting of the
elements of $\Sigma (\Phi _{\lambda })$ which preserves node labels, is
isomorphic to $\Sigma _{e_{1}}\times \ldots \times \Sigma _{e_{p}}$. It
induces all isomorphisms among these chequered diagrams. Of course $|\mathrm{Aut}(D)|=\prod_{1\leq i\leq p}e_{i}!$. We have almost established the first step
in the proof of Theorem \ref{Theorem_Z2_boson}, namely

\begin{proposition}
\label{Prop_Z2boson_cheq} We have 
\begin{equation}
Z_{M}^{(2)}(\tau _{1},\tau _{2},\epsilon )=Z_{M}^{(1)}(\tau
_{1})Z_{M}^{(1)}(\tau _{2})\sum_{D}\frac{\gamma (D)}{|\mathrm{Aut}(D)|},
\label{eq: part func D}
\end{equation}%
where $D$ ranges over isomorphism classes of chequered configurations and 
\begin{equation}
\gamma (D)=\frac{E(\lambda )}{\prod_{i}i^{e_{i}}}\epsilon ^{\sum ie_{i}}.
\label{gammaD}
\end{equation}
\end{proposition}

Proposition \ref{Prop_Z2boson_cheq} follows from what we have said together
with (\ref{eq: part func lambda}). It is only necessary to point out that
because the label subgroup induces all isomorphisms of chequered diagrams,
when we sum over isomorphism classes of such diagrams in (\ref{eq: part func
lambda}) the term $\prod_{i}e_{i}!$ must be replaced by $|\mathrm{Aut}(D)|$. 
$\square $

\bigskip Recalling the weights (\ref{eq: eqomega}), we define 
\begin{equation*}
\zeta (D)=\Pi _{E}\zeta (E),
\end{equation*}%
where the product is taken over the edges $E$ of $D$ and $\zeta(E)$ is as in
(\ref{omegaval}).

\begin{lemma}
\label{Lemma_om_gamma} For all $D$ we have%
\begin{equation}
\zeta (D)=\gamma (D).  \label{omgamma}
\end{equation}
\end{lemma}

\noindent \textbf{Proof.} Let $D$ be determined by a partition $\lambda
=\{1^{e_{1}}\ldots p^{e_{p}}\}$ and a pair of involutions $\phi ,\psi \in
F(\Phi _{\lambda })$, and let $(a,b),(r,s)$ range over the orbits of $\phi $
resp. $\psi $ on $\Phi _{\lambda }$. Then we find 
\begin{eqnarray*}
\frac{E(\lambda )}{\prod_{i}i^{e_{i}}}\epsilon ^{\sum ie_{i}} &=&\frac{%
\prod_{(a,b)}C(a,b,\tau _{1})\prod_{(r,s)}C(r,s,\tau _{2})}{%
\prod_{i}i^{e_{i}}}\epsilon ^{\sum ie_{i}} \\
&=&\prod_{(ab)}\frac{\epsilon ^{(a+b)/2}}{\sqrt{ab}}C(a,b,\tau
_{1})\prod_{(rs)}\frac{\epsilon ^{(r+s)/2}}{\sqrt{rs}}C(r,s,\tau _{2}) \\
&=&\prod_{(ab)}A_{1}(a,b)\prod_{(rs)}A_{2}(r,s)=\zeta (D).\ \ \ \square
\end{eqnarray*}

We may represent a chequered diagram formally as a product 
\begin{equation}
D=\prod_{i}L_{i}^{m_{i}}  \label{DLprod}
\end{equation}%
in case $D$ is the disjoint union of unoriented chequered cycles $L_{i}$
with multiplicity $m_{i}$. Then $\mathrm{Aut}(D)$ is isomorphic to the
direct product of the groups 
$\mathrm{Aut}(L_{i}^{m_{i}})$
of order $\left\vert \mathrm{Aut}(L_{i}^{m_{i}})\right\vert =\left\vert 
\mathrm{Aut}(L_{i})\right\vert ^{m_{i}}m_{i}!$ 
so that 
\begin{equation*}
|\mathrm{Aut}(D)|=\prod_{i}|\mathrm{Aut}(L_{i}^{m_{i}})|m_i!.
\end{equation*}%
Noting that the expression $\zeta (D)$ is multiplicative over disjoint
unions of diagrams, we calculate

\begin{eqnarray*}
\sum_{D}\frac{\zeta (D)}{|\mathrm{Aut}(D)|} 
%&=&\prod_{L}\sum_{k\geq 0}\frac{\zeta (L^{k})}{|\mathrm{Aut}(L^{k})|}\\
&=&\prod_{L}\sum_{k\geq 0}\frac{\zeta (L)^{k}}{|\mathrm{Aut}(L)|^{k}k!} \\
&=&\prod_{L}\exp \left( \frac{\zeta (L)}{|\mathrm{Aut}(L)|}\right) \\
&=&\exp \left( \sum_{L}\frac{\zeta (L)}{|\mathrm{Aut}(L)|}\right) ,
\end{eqnarray*}%
where $L$ ranges over isomorphism classes of unoriented chequered cycles.
Now $\mathrm{Aut}(L)$ is either a dihedral group of order $2r$ or a cyclic
group of order $r$ for some $r\geq 1$, depending on whether $L$ admits a
reflection symmetry or not. If we now \emph{orient} our cycles, say in a
clockwise direction, then we can replace the previous sum over $L$ by a sum
over the set of (isomorphism classes of) \emph{oriented} chequered cycles $%
\mathcal{O}$ to obtain 
\begin{equation}
\sum_{D}\frac{\zeta (D)}{|\mathrm{Aut}(D)|}=\exp \left( \frac{1}{2}%
\sum_{M\in \mathcal{O}}\frac{\zeta (M)}{|\mathrm{Aut}(M)|}\right) .
\label{omDM}
\end{equation}

Let $\mathcal{O}_{2n}\subset $ $\mathcal{O}$ denoted the set of oriented
chequered cycles with $2n$ nodes. Then we have

\begin{lemma}
\label{Lemma_TrA1A2}%
\begin{equation}
\mathrm{Tr}((A_{1}A_{2})^{n})=\sum_{M\in \mathcal{O}_{2n}}\frac{n}{|\mathrm{%
Aut}(M)|}\zeta (M).  \label{trace}
\end{equation}
\end{lemma}

\noindent \textbf{Proof.} The contribution $%
A_{1}(i_{1},i_{2})A_{2}(i_{2},i_{3})\ldots A_{2}(i_{2n},i_{1})$ to the
left-hand-side of (\ref{trace}) is equal to the weight $\zeta (M)$ for some $%
M\in \mathcal{O}_{2n}$ with vertices $i_{1},i_{2},\ldots i_{2n}$. Let $%
\sigma =\left( 
\begin{array}{ccccc}
i_{1} & \ldots & i_{k} & \ldots & i_{2n} \\ 
i_{3} & \ldots & i_{k+2} & \ldots & i_{2}%
\end{array}%
\right) $ denote the order $n$ permutation of the indices which generates
rotations of $M$. Then $\mathrm{Aut}(M)=\langle \sigma ^{m}\rangle $ for
some $m=n/|\mathrm{Aut}(M)|$. Now sum over all $i_{k}$ to compute $\mathrm{Tr%
}((A_{1}A_{2})^{n})$, noting that for inequivalent $M$ the weight $\zeta (M)$
occurs with multiplicity $m$. The Lemma follows. \hfill $\square $

\medskip We may now complete the proof of Theorem \ref{Theorem_Z2_boson}.
From (\ref{omDM}) and (\ref{trace}) we obtain 
\begin{eqnarray*}
\sum_{D}\frac{\zeta (D)}{|\mathrm{Aut}(D)|} &=&\exp \left( \frac{1}{2}%
\mathrm{Tr}(\sum_{n}\frac{1}{n}(A_{1}A_{2})^{n})\right) \\
&=&\exp (-\frac{1}{2}\mathrm{Tr}(\log (1-A_{1}A_{2}))) \\
&=&\det (\exp (-\frac{1}{2}(\log (1-A_{1}A_{2})))) \\
&=&(\det (1-A_{1}A_{2}))^{-1/2}.\quad \square
\end{eqnarray*}

\bigskip We may also obtain a product formula for $Z_{M}^{(2)}(\tau
_{1},\tau _{2},\epsilon )$ as follows. Recalling the notation (\ref%
{mathcalRdef}), for each oriented chequered cycle $M$, $\mathrm{Aut}(M)$ is
a cyclic group of order $r$ for some $r\geq 1$. Furthermore it is evident
that there is a rotationless chequered cycle $N$ with $\zeta (M)=\zeta
(N)^{r}$. Indeed, $N$ may be obtained by taking a suitable consecutive
sequence of $n/r$ nodes of $M$, where $n$ is the total number of nodes of $M$%
. We thus see that 
\begin{eqnarray*}
\sum_{M\in \mathcal{O}}\frac{\zeta (M)}{|\mathrm{Aut}(M)|} &=&\sum_{N\in 
\mathcal{R}}\sum_{r\geq 1}\frac{\zeta (N)^{r}}{r} \\
&=&-\sum_{N\in \mathcal{R}}\log (1-\zeta (N)).
\end{eqnarray*}%
Then (\ref{omDM}) implies%
\begin{equation}
\det (1-A_{1}A_{2})=\prod_{N\in \mathcal{R}}(1-\zeta (N)),  \label{detprod}
\end{equation}%
\ and thus we obtain

\begin{theorem}
\label{Theorem_Z2_boson_prod} Let $M$ be the vertex operator algebra of one
free boson. Then 
\begin{equation}
Z_{M}^{(2)}(\tau _{1},\tau _{2},\epsilon )= \frac{Z_{M}^{(1)}(\tau
_{1})Z_{M}^{(1)}(\tau _{2})}{\prod_{N\in \mathcal{R}}(1-\zeta (N))^{1/2}}.
\label{eq: part func free bos}
\end{equation}
\end{theorem}

\subsection{Holomorphic and Modular Invariance Properties}

\label{subsect_modular}

In Section \ref{subsect_epsilon} we reviewed the genus two $\epsilon$-sewing
formalism and introduced the domain $\mathcal{D}^{\epsilon }$ parameterizing
the genus two surface. An immediate consequence of Theorem \ref%
{Theorem_Z2_boson} and Theorem \ref{Theorem_A1A2}(b) is the following:

\begin{theorem}
\label{Theorem_Z2_boson_eps_hol} $Z_{M}^{(2)}(\tau _{1},\tau _{2},\epsilon )$
is holomorphic on\ the domain $\mathcal{D}^{\epsilon }$. \hfill $\square $
\end{theorem}

We next consider the automorphic properties of the genus two partition
function with respect to the group $G$ reviewed in Section~\ref%
{subsect_epsilon}. For two free bosons the genus one partition function is 
\begin{equation}
Z_{M^{2}}^{(1)}(\tau )=\frac{1}{\eta (\tau )^{2}}.  \label{etafunceqn}
\end{equation}%
Let $\chi $ be the character of $SL(2,\mathbb{Z})$ defined by its action on $%
\eta (\tau )^{-2}$, i.e. 
\begin{equation}
\eta (\gamma \tau )^{-2}=\chi (\gamma )\eta (\tau )^{-2}(c\tau +d)^{-1},
\label{eq: eta}
\end{equation}%
where $\gamma =\left( 
\begin{array}{cc}
a & b \\ 
c & d%
\end{array}%
\right) \in SL(2,\mathbb{Z})$. Recall (e.g. \cite{Se}) that $\chi (\gamma )$
is a twelfth root of unity. For a function $f(\tau )$ on $\mathbb{H}%
_{1},k\in \mathbb{Z}$ and $\gamma \in SL(2,\mathbb{Z})$, we define 
\begin{equation}
f(\tau )|_{k}\gamma =f(\gamma \tau )\ (c\tau +d)^{-k},  \label{slashaction}
\end{equation}%
so that 
\begin{equation}
Z_{M^{2}}^{(1)}(\tau )|_{-1}\gamma =\chi (\gamma )Z_{M^{2}}^{(1)}(\tau ).
\label{Z1modgam}
\end{equation}

The genus two partition function for two free bosons is 
\begin{equation}
Z_{M^{2}}^{(2)}(\tau _{1},\tau _{2},\epsilon ) =\frac{1}{\eta (\tau
_{1})^{2}\eta (\tau _{2})^{2}\det (I-A_{1}A_{2})}.  \label{Z2mod_eps}
\end{equation}%
Analogously to (\ref{slashaction}), we define 
\begin{equation}
f(\tau _{1},\tau _{2},\epsilon )|_{k}\gamma =f(\gamma (\tau _{1},\tau
_{2},\epsilon ))\det (C\Omega +D)^{-k}.  \label{eq: Gaction}
\end{equation}%
Here, the action of $\gamma $ on the right-hand-side is as in (\ref{GDeps}).
We have abused notation by adopting the following conventions in (\ref{eq:
Gaction}), which we continue to use below: 
\begin{equation}
\Omega =F^{\epsilon }(\tau _{1},\tau _{2},\epsilon ),\ \gamma =\left( 
\begin{array}{cc}
A & B \\ 
C & D%
\end{array}%
\right) \in Sp(4,\mathbb{Z})  \label{Omegaconvention}
\end{equation}%
where $F^{\epsilon }$ is as in Theorem \ref{Theorem_period_eps}, and $\gamma 
$ is identified with an element of $Sp(4,\mathbb{Z})$ via (\ref{Gamma1Gamma2}%
)-(\ref{betagen}). Note that (\ref{eq: Gaction}) defines a right action of $%
G $ on functions $f(\tau _{1},\tau _{2},\epsilon )$. We will establish the
natural extension of (\ref{Z1modgam}) to the genus $2$ case. To describe this, introduce the
character $\chi ^{(2)}$ of $G$ defined by 
\begin{equation*}
\chi ^{(2)}(\gamma _{1}\gamma _{2}\beta ^{m})=(-1)^{m}\chi (\gamma
_{1}\gamma _{2}),\quad \ \gamma _{i}\in \Gamma _{i},\ i=1,2,
\end{equation*}%
(notation as in (\ref{Gamma1Gamma2}), (\ref{betagen})). Thus $\chi ^{(2)}$ takes
values which are twelfth roots of unity, and we have

\begin{theorem}
\label{Theorem_Z2_G}If $\gamma \in G$ then 
\begin{equation*}
Z_{M^{2}}^{(2)}(\tau _{1},\tau _{2},\epsilon )|_{-1}\gamma =\chi
^{(2)}(\gamma )Z_{M^{2}}^{(2)}(\tau _{1},\tau _{2},\epsilon ).
\end{equation*}
\end{theorem}

\begin{corollary}
\label{Cor_Z2_24G} For the rank 24 Heisenberg VOA $M^{24}$ we have 
\begin{equation*}
Z_{M^{24}}^{(2)}(\tau _{1},\tau _{2},\epsilon )|_{-12}\gamma
=Z_{M^{24}}^{(2)}(\tau _{1},\tau _{2},\epsilon ),
\end{equation*}
for $\gamma \in G$.
\end{corollary}

\noindent \textbf{Proof.} We will give two different proofs of this result.
Using the convention (\ref{Omegaconvention}), we have to show that 
\begin{equation}
Z_{M^{2}}^{(2)}(\gamma (\tau _{1},\tau _{2},\epsilon ))\det (C\Omega
+D)=\chi ^{(2)}(\gamma )Z_{M^{2}}^{(2)}(\tau _{1},\tau _{2},\epsilon )
\label{eq: Z^2identity}
\end{equation}%
for $\gamma \in G$, and it is enough to do this for a generating set of $G$.
If $\gamma =\beta $ then the result is clear since $\det (C\Omega +D)=\chi
^{(2)}(\beta )=-1$ and $\beta $ exchanges $\tau _{1}$ and $\tau _{2}$. So we
may assume that $\gamma =(\gamma _{1},\gamma _{2})\in \Gamma _{1}$ $\times $ 
$\Gamma _{2}$.

\medskip Our first proof utilizes the determinant formula (\ref{Z2_1bos}) as
follows. For $\gamma _{1}\in \Gamma _{1}$, define $A_{a}^{\prime }(k,l,\tau
_{a},\epsilon )=A_{a}(k,l,\gamma _{1}\tau _{a},\frac{\epsilon }{c_{1}\tau
_{1}+d_{1}})$ following (\ref{GDeps}). We find from Section~4.4 of \cite{MT2}
that 
\begin{eqnarray*}
I-A_{1}^{\prime }A_{2}^{\prime } &=&I-A_{1}A_{2}-\kappa \Delta A_{2} \\
&=&(I-\kappa S).(I-A_{1}A_{2}),
\end{eqnarray*}%
where $\Delta (k,l) =\delta _{k1}\delta _{l1}$, $\kappa -\frac{\epsilon }{%
2\pi i}\frac{c_{1}}{c_{1}\tau _{1}+d_{1}}$ and $S(k,l)=\delta
_{k1}(A_{2}(I-A_{1}A_{2})^{-1})(1,l)$. Since $\det (I-A_{1}A_{2})$ and $\det
(I-A_{1}^{\prime }A_{2}^{\prime })$ are convergent on $\mathcal{D}^{\epsilon
}$ we find%
\begin{equation*}
\det (I-A_{1}^{\prime }A_{2}^{\prime })=\det (I-\kappa S)\det (I-A_{1}A_{2}).
\end{equation*}%
But $\det (I-\kappa S)=1-\kappa S(1,1)=\frac{c_{1}\Omega _{11}+d_{1}}{
c_{1}\tau _{1}+d_{1}}$ which implies (\ref{eq: Z^2identity}) for $\gamma
_{1}\in \Gamma _{1}$. A similar proof applies for $\gamma _{2}\in \Gamma
_{2} $.

\medskip The second proof uses Proposition \ref{Prop_Om12_R21expansion}
together with (\ref{eq: omega_12prod}), which tell us that 
\begin{equation}
Z_{M^{2}}^{(2)}(\tau _{1},\tau _{2},\epsilon )=\frac{-2\pi i\Omega _{12}}{%
\epsilon \eta (\tau _{1})^{2}\eta (\tau _{2})^{2}}\prod_{\mathcal{R}^{\prime
}}(1-\zeta (L))^{-1},  \label{eq: Z^2 omega_12}
\end{equation}%
where $\mathcal{R}^{\prime }=\mathcal{R}\setminus \mathcal{R}_{21}$. Now in
general a term $\zeta (L)$ will not be invariant under the action of $\gamma 
$. This is because of the presence of quasi-modular terms $A_{a}(1,1)$
arising from $E_{2}(\tau _{a})$. But it is clear from (\ref{GDeps}) and the
definition (\ref{Ckldef}) of $C(k,l,\tau )$ together with its
modular-invariance properties that if $L\in \mathcal{R}^{\prime }$ then such
terms are absent and $\zeta (L)$ \emph{is} invariant. So the product term in
(\ref{eq: Z^2 omega_12}) is invariant under the action of $\gamma $.

Next, we see from (\ref{GDeps}) that the expression $\epsilon \eta (\tau
_{1})^{2}\eta (\tau _{2})^{2}$ is invariant under the action of $\gamma $ up
to a scalar $\chi (\gamma _{1})\chi (\gamma _{2})=\chi ^{(2)}(\gamma )$.
This reduces the proof of (\ref{eq: Z^2identity}) to showing that 
\begin{equation*}
(\gamma _{1},\gamma _{2}):\Omega _{12}\mapsto \Omega _{12}\ \det (C\Omega
+D)^{-1},
\end{equation*}%
and this is implicit in (\ref{eq: modtrans}) upon applying Theorem \ref%
{TheoremGequiv}. This completes the second proof of Theorem \ref%
{Theorem_Z2_G}. \hfill $\square$

\begin{remark}
\label{automorphyremark} An unusual feature of the formulas in Theorem \ref%
{Theorem_Z2_G} and Corollary \ref{Cor_Z2_24G} is that the definition of the
automorphy factor $\det (C\Omega +D)$ requires the map $F^{\epsilon }:%
\mathcal{D}^{\epsilon }\rightarrow \mathbb{H}_{2}$.
\end{remark}

\subsection{Some Genus Two $n$-Point Functions}

\label{subsect_Heisenberg_npt} In this section we calculate some examples of
genus two $n$-point functions for the rank one Heisenberg VOA $M$. A general
analysis of all such functions will appear elsewhere \cite{MT5}. We consider
here the examples of the $n$-point function for the Heisenberg vector $a$
and the 1-point function for the Virasoro vector $\tilde\omega$. We find
that the formal differential form (\ref{F2}) associated with the Heisenberg $%
n$-point function is described in terms of the global symmetric two form $%
\omega^{(2)}$ \cite{TUY} whereas the Virasoro 1-point function is described
by the genus two projective connection \cite{Gu}. These results illustrate
the general conjecture made in Remark \ref{rem:F2global}.

We first consider the example of the Heisenberg vector 1-point function
where $a$ is inserted at $x$ on the left torus (say). Since $%
Z_{M}^{(1)}(Y[a,x]v;\tau)=0$ for a Fock vector $v=v(\lambda)$ for even $%
\vert\lambda\vert$ and $Z_{M}^{(1)}(v;\tau)=0$ for odd $\vert\lambda\vert$ 
\cite{MT1} we find from (\ref{Z2n_pt_eps}) that $Z_{M}^{(2)}(a,x\vert\tau_1,%
\tau_2,\epsilon)=0$.

Consider next the 2-point function for two Heisenberg vectors inserted on
the left torus at $x_1,x_2\in \hat{\mathcal{S}_1}$ with 
\begin{equation}
Z_{M}^{(2)}(a,x_1;a,x_2\vert \tau_1,\tau_2,\epsilon)= \sum_{r\geq
0}\epsilon^{r}\sum_{v\in M_{[r]}}Z_{M}^{(1)}(Y[a,x_1]Y[a,x_2]v;\tau
_{1})Z_{M}^{(1)}(\bar{v};\tau _{2}),  \label{Z2a1a1}
\end{equation}%
Following (\ref{F2}) of Proposition \ref{Prop_Z2sew}, we consider the
associated formal differential form $\mathcal{F}^{(2)}(a,a;\tau_1,\tau_2,%
\epsilon)$ for (\ref{Z2a1a1}) and find that it is determined by the bilinear
form $\omega^{(2)}$ of (\ref{omegag}):

\begin{theorem}
\label{theorem:Z2aa} The genus two Heisenberg vector 2-point function is 
\begin{eqnarray}
\mathcal{F}_{M}^{(2)}(a,a;\tau_1,\tau_2,\epsilon)=\omega^{(2)}
Z_{M}^{(2)}(\tau_1,\tau_2,\epsilon).  \label{eq: F2aa}
\end{eqnarray}
\end{theorem}

\noindent \textbf{Proof.} The proof proceeds along the same lines as Theorem %
\ref{Theorem_Z2_boson}. As before, we let $v(\lambda )$ denote a Heisenberg
Fock vector (\ref{eq: sq fock vec}) determined by an unrestricted partition $%
\lambda =\{1^{e_{1}}\ldots p^{e_{p}}\}$ with label set $\Phi _{\lambda }$.
Define a label set for the three vectors $a,a,v(\lambda )$ given by $%
\Phi=\Phi_1\cup \Phi_2\cup \Phi_3$ for $\Phi_1,\Phi_2=\{ 1\}$ and $%
\Phi_3=\Phi_{\lambda }$ and let $F(\Phi)$ denote the set of fixed point free
involutions on $\Phi$. For $\phi=\ldots (r s)\ldots\in F(\Phi)$ let $%
\Gamma_1 (x_1,x_2,\phi)=\prod_{(r,s)}\gamma(r,s)$ where for $r\in \Phi_i$
and $s\in \Phi_j$ 
\begin{equation}
\gamma (r,s)=\left \{ 
\begin{array}{ll}
D(1,1,x_1-x_2,\tau_1 )=P_2(\tau_1,x_1-x_2), & i=1;\ j=2 \\ 
D(1,s,x_{i},\tau_1 )=sP_{s+1}(\tau_1,x_i), & i=1,2;\ j=3 \\ 
C(r,s,\tau_1 ), & i,j=3,%
\end{array}
\right.  \label{eq: X12 weight}
\end{equation}
for $C,D$ of (\ref{Ckldef}) and (\ref{Dkldef}). Then following Corollary 1
of \cite{MT1} we find for even $\vert\lambda\vert$ that 
\begin{equation*}
Z_{M}^{(1)}(Y[a,x_1]Y[a,x_2]v(\lambda ),\tau_1 )=Z_{M}^{(1)}(\tau_1 )
\sum_{\phi \in F(\Phi)}\Gamma_1 (x_1,x_2,\phi ).
\end{equation*}%
Recalling that $\mathcal{F}^{(2)}(a,a;\tau_1,\tau_2,%
\epsilon)=Z_{M}^{(2)}(a,x_1;a,x_2\vert\tau _{1},\tau _{2},\epsilon )dx_1dx_2$
we then obtain the following analogue of (\ref{eq: part func lambda}) 
\begin{equation}
\mathcal{F}^{(2)}(a,a;\tau_1,\tau_2,\epsilon)=
Z_{M}^{(1)}(\tau_{1})Z_{M}^{(1)}(\tau _{2})\sum_{\lambda =\{i^{e_{i}}\}}%
\frac{E(x_1,x_2,\lambda)}{\prod_{i}i^{e_{i}}e_{i}!}\epsilon^{\sum
ie_{i}}dx_1dx_2,  \label{eq: part func lambda a1a2}
\end{equation}%
where 
\begin{equation*}
E(x_1,x_2,\lambda ) = \sum_{\phi\in F(\Phi),\ \psi \in F(\Phi _{\lambda
})}\Gamma _{1}(x_1,x_2,\phi)\Gamma _{2}(\psi ),
\end{equation*}
with $\Gamma _{2}(\psi)$ as before.

The expression (\ref{eq: part func lambda a1a2}) can be interpreted as a sum
of weights $\zeta (D)$ associated with isomorphism classes of chequered
configurations $D$ where, in this case, \textbf{each} configuration includes
two distinguished valence one nodes of type $1,x_i$ (see Section~\ref%
{subsect_Omega_graph}) corresponding to the label sets $\Phi_1,\Phi_2=\{ 1\}$%
. As before, $\zeta (D)=\prod_E\zeta(E)$ for standard chequered edges $E$ (%
\ref{omegaval}) augmented by the contributions for edges connected to the
two valence one nodes with weights as in (\ref{eq: zeta1a2}) (for $a=1$).
Then we find, as in Proposition \ref{Prop_Z2boson_cheq}, that 
\begin{equation*}
\mathcal{F}^{(2)}(a,a;\tau_1,\tau_2,\epsilon)=Z_{M}^{(1)}(\tau
_{1})Z_{M}^{(1)}(\tau _{2})\sum_{D}\frac{\zeta (D)}{\prod_{i}e_{i}!}dx_1dx_2.
\end{equation*}
Each $D$ can be decomposed into \emph{exactly} one necklace configuration $N$
of type $\mathcal{N}_{11}^{x_1,x_2}$ of (\ref{Nabxydef}) connecting the two
distinguished nodes and a standard configuration ${\hat D}$ of the type
appearing in Subsection~\ref{subsect_BosonZ2} so that $\zeta (D)=\zeta
(N)\zeta ({\hat D})$. Furthermore, if $\lambda^\prime=\{1^{e^\prime_{1}}%
\ldots p^{e^\prime_{p}}\}$ is the subset of $\lambda$ that labels ${\hat D}$
then the necklace contribution $\zeta (N)$ occurs with multiplicity $%
\prod_{i}\frac{e_{i}!}{e^\prime_{i}!}=\frac{|\mathrm{Aut}({\ D})|}{|\mathrm{%
Aut}({\hat D})|}$. It follows that 
\begin{eqnarray*}
\mathcal{F}^{(2)}(a,a;\tau_1,\tau_2,\epsilon)&=&Z_{M}^{(1)}(\tau
_{1})Z_{M}^{(1)}(\tau _{2}) \sum_{{\hat D}}\frac{\zeta ({\hat D})}{|\mathrm{%
Aut}({\hat D})|} \sum_{N\in \mathcal{N}_{11}^{x_1,x_2}}\zeta(N)dx_1dx_2 \\
&=&Z_{M}^{(2)}(\tau _{1},\tau _{2},\epsilon )\zeta_{11}^{x_1,x_2}dx_1dx_2 \\
&=&Z_{M}^{(2)}(\tau _{1},\tau _{2},\epsilon )\omega^{(2)}(x_1,x_2),
\end{eqnarray*}
using (\ref{om_graph}) of Proposition \ref{Propomeganugraph}. Applying
Proposition \ref{Prop_Z2sew}, the same two form arises for the other
possible insertions of two Heisenberg vectors. Alternatively, a similar
explicit calculation can be carried out in each case leading to the
expressions for $\omega^{(2)}$ described by (\ref{om_graph}). \hfill $%
\square $

In a similar fashion one can generally show that the $n$-point function for $%
n$ Heisenberg vectors vanishes for $n$ odd and for $n$ even is determined by
the global symmetric meromorphic $n$ form given by the symmetric (tensor)
product 
\begin{equation}
\mathrm{Sym}_n\omega^{(2)}= \sum_{\psi}\prod_{(r,s)}\omega^{(2)}(x_r,x_s),
\label{eq: bos n form}
\end{equation}
where the sum is taken over the set of fixed point free involutions $%
\psi=\ldots (rs)\ldots$ of the labels $\{1,\ldots,n\}$. Then one finds

\begin{theorem}
\label{theorem:Z2an} The genus two Heisenberg vector $n$-point function is
given by the global symmetric meromorphic $n$-form 
\begin{equation}
\mathcal{F}_{M}^{(2)}(a,\ldots,a;\tau_1,\tau_2,\epsilon)=\mathrm{Sym}%
_n\omega^{(2)} Z_{M}^{(2)}(\tau_1,\tau_2,\epsilon).  \label{eq: F2n}
\end{equation}
\end{theorem}

Theorem \ref{theorem:Z2an} is in agreement with earlier results in \cite{TUY}
based on an assumed analytic structure for the ratio $\mathcal{F}%
_{M}^{(2)}(a,\ldots,a;\tau_1,\tau_2,\epsilon)/Z_{M}^{(2)}(\tau_1,\tau_2,%
\epsilon)$.

Using the associativity property of a VOA, the genus two Heisenberg $n$%
-point function (\ref{eq: F2n}) is a generator of \emph{all} genus two $n$%
-point functions for $M$ in an analagous way to that described for genus one
in  \cite{MT1}. This will be further developed elsewhere \cite{MT5}. We
illustrate this by computing the 1-point function for the Virasoro vector ${%
\tilde \omega}=\frac{1}{2}a[-1]a$. This is determined by the genus two
projective connection defined by e.g. \cite{Gu} 
\begin{equation}
s^{(2)}(x)=6\lim_{x \rightarrow y} \left (\omega^{(2)}(x,y)-\frac{dxdy}{%
(x-y)^2}\right ).  \label{eq:proj_con}
\end{equation}
We then find

\begin{proposition}
\label{Prop:Z2omega} The genus two 1-point function for the Virasoro vector $%
\tilde\omega$ is 
\begin{eqnarray}
\mathcal{F}_{M}^{(2)}({\tilde\omega};\tau_1,\tau_2,\epsilon)= \frac{1}{12}%
s^{(2)}Z_{M}^{(2)}(\tau_1,\tau_2,\epsilon).  \label{eq: F2omega}
\end{eqnarray}
\end{proposition}

\noindent \textbf{Proof.} Using the associativity property of a VOA we have 
\cite{MT1} 
\begin{eqnarray*}
Z_{M}^{(1)}(Y[a,x_1]Y[a,x_2]v;\tau_{1})&=&Z_{M}^{(1)}(Y[Y[a,x_1-x_2]a,x_2]v;%
\tau_{1}) \\
&=&\frac{Z_{M}^{(1)}(v;\tau_{1})}{(x_1-x_2)^{2}} +2Z_{M}^{(1)}(Y[{\tilde
\omega},x_2]v;\tau_{1})+\ldots
\end{eqnarray*}
Hence using the Heisenberg 2-point function (\ref{Z2a1a1}) we find 
\begin{eqnarray*}
\mathcal{F}^{(2)}(\omega;\tau_1,\tau_2,\epsilon)&=& \lim_{x_1\rightarrow x_2}%
\frac{1}{2}\left (Z_{M}^{(2)}(a,x_1;a,x_2\vert \tau_1,\tau_2,\epsilon) -%
\frac{Z_{M}^{(2)}(\tau_1,\tau_2,\epsilon)}{(x_1-x_2)^2}\right)dx_1dx_2 \\
&=&\frac{1}{12}s^{(2)}(x_1)Z_{M}^{(2)}(\tau_1,\tau_2,\epsilon).\ \ \ \square
\end{eqnarray*}

Notice that $\mathcal{F}^{(2)}(\omega;\tau_1,\tau_2,\epsilon)$ is not a
global differential 2-form since $s^{(2)}(x)$ transforms under a general
conformal transformation $\phi(x)$ (\cite{Gu}) as 
\begin{equation}
s^{(2)}(\phi(x))=s^{(2)}(x)-\{\phi;x\}dx^2,  \label{eq: proj_con_fx}
\end{equation}
where $\{\phi;x\}=\frac{\phi^{\prime\prime\prime}}{%
\phi^{\prime}}-\frac{3}{2}\left(\frac{\phi^{\prime\prime}}{\phi^{\prime}}%
\right)^2$ is the usual Schwarzian derivative. This property of the Virasoro 1-point function has previously
been discussed many times in the physics and mathematics literature based on
a variety of stronger assumptions e.g. \cite{EO, TUY, FS, U, Z2}.

%%%%%%%% end n pt fun

\section{Heisenberg Modules, Lattice VOAs and Theta Series}

\label{sect_Lattice} In this Section we generalize the methods of Section~%
\ref{sect_Boson} to compute the genus two partition function for a pair of
Heisenberg modules. We consider the genus two $n$-point function for the
Heisenberg vector and the Virasoro 1-point function. We apply these result
to obtain closed formulas for the genus two partition function for a lattice
VOA $V_L$ (in terms of the genus two Siegel theta function for $L$) and the
\lq twisted\rq\ genus two partition function for the 
$\mathbb Z$
-lattice VOA (in terms of the genus two character valued Riemann theta
function). We finally derive a genus two Ward identity for the Virasoro
1-point function for these theories.

\subsection{Heisenberg Modules}

\label{subsect_Hmodules} In this section we discuss the genus two partition
function for a pair of simple Heisenberg modules $M\otimes e^{\alpha_1}$ and 
$M\otimes e^{\alpha_2}$ for $\alpha_1,\alpha_2\in {}$. The partition
function is then 
\begin{equation}
Z_{\alpha_1,\alpha_2}^{(2)}(\tau _{1},\tau _{2},\epsilon )=\sum_{n\geq
0}\epsilon ^{n}\sum_{u\in M_{[n]}}Z_{M\otimes e^{\alpha_1 }}^{(1)}(u,\tau
_{1})Z_{M\otimes e^{\alpha_2 }}^{(1)}({\bar u},\tau _{2}),
\label{eq: alpha part func}
\end{equation}
where $u$ ranges over any basis for $M_{[n]}$. An explicit formula for $%
Z_{M\otimes e^{\alpha}}^{(1)}(u,\tau )$ was given in \cite{MT2} (Corollary 3
and Theorem 1). We are going to use these results, together with graphical
techniques similar to those employed for free bosons in Section~\ref%
{sect_Boson} to establish a closed formula for (\ref{eq: alpha part func}).
Letting ${\alpha}. \Omega.{\alpha}=\sum_{i,j=1,2}{\alpha_i} \Omega_{ij}{%
\alpha_j}$ where $\Omega_{ij}$ is the genus two period matrix we find

\begin{theorem}
\label{Theorem_Z2_Malpha} We have 
\begin{equation}
Z_{\alpha_1,\alpha_2}^{(2)}(\tau _{1},\tau _{2},\epsilon )=e^{i\pi
\alpha.\Omega.\alpha}Z_{M}^{(2)}(\tau_{1},\tau _{2},\epsilon ).
\label{eq: Z2alphaOmega}
\end{equation}
$Z_{\alpha_1,\alpha_2}^{(2)}(\tau _{1},\tau _{2},\epsilon )$ is holomorphic
on the domain $\mathcal{D}^{\epsilon}$.
\end{theorem}

\begin{remark}
This is a natural generalization of the genus one partition function
relation $Z_{M\otimes
e^{\alpha}}^{(1)}(\tau)=q^{\alpha^2/2}Z_{M}^{(1)}(\tau) $.
\end{remark}

\noindent \textbf{Proof.} Consider the Fock basis vectors $v=v(\lambda)$
(cf. (\ref{eq: sq fock vec})) identified with partitions $\lambda
=\{i^{e_{i}}\}$ as in Section~\ref{sect_Boson}. Recall that $\lambda $
defines a labeled set $\Phi _{\lambda }$ with $e_{i}$ nodes labeled $i$. It
is useful to re-state Corollary 3 of \cite{MT1} in the following form: 
\begin{equation}
Z_{M\otimes e^{\alpha}}^{(1)}(v,\tau )=Z_{M}^{(1)}(\tau
)q^{\alpha^2/2}\sum_{\phi }\Gamma _{\lambda ,\alpha }(\phi ).
\label{eq: Z^1Malpha}
\end{equation}%
Here, $\phi $ ranges over the set of involutions 
\begin{equation}
\mathrm{Inv}_{1}(\Phi _{\lambda })=\{\phi \in \mathrm{Inv}(\Phi _{\lambda
})|\ p\in \mathrm{Fix}(\phi )\Rightarrow p\mbox{ has label}\ 1\}.
\label{eq: inv}
\end{equation}%
In words, $\phi $ is an involution in the symmetric group $\Sigma (\Phi
_{\lambda })$ such that all fixed-points of $\phi $ carry the label $1$.
Note that this includes the fixed-point-free involutions, which were the
only involutions which played a role in the case of free bosons. The main
difference between the free bosonic VOA and its modules is the need to
include additional involutions in the latter case. In particular, we note
that permutations with $\vert \lambda\vert$ odd can contribute in this case
for $\lambda =\{i^{e_{i}}\}$ with $e_1$ odd. Finally, 
\begin{equation}
\Gamma _{\lambda ,\alpha }(\phi ,\tau )=\Gamma _{\lambda ,\alpha }(\phi
)=\prod_{\Xi }\Gamma (\Xi ),  \label{eq: gammavalpha}
\end{equation}%
where $\Xi $ ranges over the orbits (of length $\leq 2$) of $\phi $ acting
on $\Phi _{\lambda }$ and 
\begin{equation}
\Gamma (\Xi )=\left\{ 
\begin{array}{ll}
C(r,s,\tau ), & \mbox{if $\Xi = \{r, s \}$}, \\ 
\alpha, & \mbox{if $\Xi = \{1 \}$}.%
\end{array}%
\right.  \label{gamma_xi}
\end{equation}
From (\ref{eq: lattice part func diag})-(\ref{eq: gammavalpha}) we get 
\begin{gather}
Z_{\alpha_1,\alpha_2}^{(2)}(\tau _{1},\tau _{2},\epsilon )=  \notag \\
Z_{M}^{(1)}(\tau _{1})Z_{M}^{(1)}(\tau _{2}) \sum_{\lambda =\{i^{e_{i}}\}} 
\frac{(-1)^{\vert\lambda\vert}E_{\alpha_1,\alpha_2 }(\lambda )} {%
\prod_{i}i^{e_{i}}e_{i}!}q_{1}^{\alpha_1^2/2}q_{2}^{\alpha_2^2/2}\epsilon
^{\sum ie_{i}},  \label{eq: newpartfunc}
\end{gather}%
where 
\begin{equation}
E_{\alpha_1,\alpha_2}(\lambda )=\sum_{\phi ,\psi \in \mathrm{Inv}_{1}(\Phi
_{\lambda })}\Gamma _{\lambda ,\alpha_1 }(\phi ,\tau _{1})\Gamma _{\lambda
,\alpha_2 }(\psi ,\tau _{2}).  \label{eq: Eeqn}
\end{equation}%
(Compare with eqns. (\ref{eq: part func lambda}) - (\ref{eq: expr E}).)

Now we follow the proof of Proposition \ref{Prop_Z2boson_cheq} to obtain an
expression analogous to (\ref{eq: part func D}), namely 
\begin{equation}
Z_{\alpha_1,\alpha_2}^{(2)}(\tau _{1},\tau _{2},\epsilon )= Z_{M}^{(1)}(\tau
_{1})Z_{M}^{(1)}(\tau _{2})\sum_{D}\frac{\gamma_{\alpha_1,\alpha_2 }^{0}(D)}{%
|\mathrm{Aut}(D)|}q_{1}^{\alpha_1^2/2}q_{2}^{\alpha_2^2/2},
\label{eq: lattice part func D}
\end{equation}%
the meaning of which we now enlarge upon. Compared to (\ref{eq: part func D}%
), the chequered diagrams $D$ which occur in (\ref{eq: lattice part func D})
are more general than before, in that they reflect the fact that the
relevant involutions may now have fixed-points. Thus $D$ is the union of its
connected (as yet unoriented) components which are either chequered cycles
as before or else chequered necklaces (see Section~\ref{subsect_Omega_graph}%
). Necklaces arise from orbits of the group $\langle \psi \phi \rangle $ on $%
\Phi _{\lambda }$ in which one of the nodes in the orbit is a fixed-point of 
$\phi $ or $\psi $. In that case the orbit will generally contain two such
nodes which comprise the end nodes of the necklace. Note that these end
nodes necessarily carry the label $1$ (cf. (\ref{eq: inv})). There is
degeneracy when \emph{both} $\phi $ and $\psi $ fix the node, in which case
the degenerate necklace obtains.

\medskip Similarly to (\ref{gammaD}), the term $\gamma_{\alpha_1 ,\alpha_2}
^{0}(D)$ in (\ref{eq: lattice part func D}) is given by 
\begin{equation}
\gamma_{\alpha_1 ,\alpha_2 }^{0}(D)=(-1)^{|\lambda|}\frac{\prod_{\Xi
_{1}}\Gamma (\Xi _{1})\prod_{\Xi _{2}}\Gamma (\Xi _{2})}{\prod_{i}i^{e_{i}}}%
\epsilon ^{\sum ie_{i}},  \label{eq: lattice gamma(D)}
\end{equation}%
where $\Xi _{1},\Xi _{2}$ range over the orbits of $\phi ,\psi $
respectively on $\Phi _{\lambda }$. As usual the summands in (\ref{eq:
lattice part func D}) are multiplicative over connected components of the
chequered diagram. This applies, in particular, to the chequered cycles
which occur, and these are independent of the lattice elements. As a result,
(\ref{eq: lattice part func D}) factors as a product of two expressions, the
first a sum over diagrams consisting only of chequered cycles and the second
a sum over diagrams consisting only of chequered necklaces. However, the
first expression corresponds precisely to the genus two partition function
for the free boson (Proposition \ref{Prop_Z2boson_cheq}). We thus obtain 
\begin{equation}
\frac{Z_{\alpha_1 ,\alpha_2}^{(2)}(\tau _{1},\tau _{2},\epsilon )}{%
Z_{M}^{(2)}(\tau _{1},\tau _{2},\epsilon )}=\sum_{D^{N}}\frac{%
\gamma_{\alpha_1 ,\alpha_2 }^{0}(D^{N})}{|\mathrm{Aut}(D^{N})|}%
q_{1}^{\alpha_1^2/2}q_{2}^{\alpha_2^2/2},  \label{eq: spliteqn}
\end{equation}%
where here $D^{N}$ ranges over all chequered diagrams all of whose connected
components are chequered necklaces. So Theorem \ref{Theorem_Z2_Malpha} is
reduced to establishing

\begin{proposition}
\label{Prop_expalpha} We have 
\begin{equation}
e^{i\pi {\alpha}. \Omega.{\alpha}}=\sum_{D^{N}}\frac{\gamma_{\alpha_1,%
\alpha_2 } ^{0}(D^{N})}{|\mathrm{Aut}(D^{N})|}q_{1}^{\alpha_1^2/2}q_{2}^{%
\alpha_2^2/2}.  \label{eq: redomega}
\end{equation}
\end{proposition}

We may apply the argument of (\ref{DLprod}) et. seq. to the inner sum in (%
\ref{eq: redomega}) to write it as an exponential expression 
\begin{equation}
\exp \{i\pi (\alpha_1^2\tau _{1}+\alpha_2^2\tau _{2})+\sum_{N}\frac{%
\gamma_{\alpha_1,\alpha_2 } ^{0}(N)}{|\mathrm{Aut}(N)|}\},
\label{eq:expgammasum}
\end{equation}%
where $N$ ranges over all unoriented chequered necklaces.

\medskip Recall the isomorphism class $\mathcal{N}_{ab}$ of oriented
chequered necklaces of type $ab$ as displayed in Fig.~3 of Section~\ref%
{subsect_Omega_graph}. Then (\ref{eq:expgammasum}) can be written as 
\begin{equation}
\exp \{i\pi (\alpha_1^2\tau _{1}+\alpha_2^2\tau _{2})+\frac{1}{2}%
\sum_{a,b\in\{1,2\}}\sum_{N_{ab}} \gamma_{\alpha_1,\alpha_2}^{0}(N_{ab})\}.
\label{eq: orexpgammasum}
\end{equation}%
where here $N$ ranges over \emph{oriented} chequered necklaces of type $ab$.

\medskip From (\ref{gamma_xi}) and (\ref{eq: lattice gamma(D)}) we see that
the contribution of the end nodes to $\gamma_{\alpha_1,\alpha_2}^{0}(N)$ is
equal to $\epsilon \alpha_{\bar a}\alpha_{\bar b}$ for a type $ab$ necklace.
The remaining edge factors of $\gamma_{\alpha_1,\alpha_2} ^{0}(N)$ have
product $\gamma(N)=\zeta (N)$ by Lemma \ref{Lemma_om_gamma}. Finally,
necklaces of type $11$ and $22$ arise from Fock vectors with an even number $%
\vert\lambda\vert$ of permutation symbols whereas necklaces of type $12$ and 
$21$ arise from Fock vectors for odd $\vert\lambda\vert$ leading to a
further $-1$ contribution in (\ref{eq: lattice gamma(D)}) in these cases.
Overall we find that 
\begin{equation*}
\sum_{N_{ab}}
\gamma_{\alpha_1,\alpha_2}^{0}(N_{ab})=(-1)^{a+b}\epsilon\alpha_{\bar
a}\alpha_{\bar b}\zeta _{ab},
\end{equation*}
recalling $\zeta _{ab}=\sum_{N\in \mathcal{N}_{ab}}\zeta (N)$. Hence (\ref%
{eq: orexpgammasum}) may be re-expressed as

\begin{equation}
\exp \{\frac{\alpha_1^2}{2}(2\pi i\tau _{1}+\epsilon \zeta _{22})+\frac{{%
\alpha_2^2}}{2}(2\pi i\tau _{2}+\epsilon \zeta _{11})-\alpha_1
\alpha_2\epsilon \zeta _{21}\},  \label{thetaLexp}
\end{equation}%
where $\zeta _{12}=\zeta _{21}$. (\ref{thetaLexp}) reproduces (\ref{eq:
g2omega}) on applying Proposition \ref{Propepsperiodgraph}.

Finally we note from Theorems \ref{Theorem_period_eps} and \ref%
{Theorem_Z2_boson_eps_hol} that $Z_{\alpha_1,\alpha_2}^{(2)}(\tau _{1},\tau
_{2},\epsilon )$ is holomorphic on the domain $\mathcal{D}^{\epsilon}$. This
completes the proof of Theorem \ref{Theorem_Z2_Malpha}. \hfill$\square$

\subsection{Some Genus Two $n$-Point Functions}

In this section we consider the genus two $n$-point functions for the
Heisenberg vector $a$ and the 1-point function for the Virasoro vector $%
\tilde\omega$ for a pair of Heisenberg modules $M\otimes e^{\alpha_i}$. We
again express each $n$-point function in terms of the associated formal
differential form following (\ref{F2}) of Proposition \ref{Prop_Z2sew}. The
results generalize those of Section~\ref{subsect_Heisenberg_npt}. They are established by
making use of similar methods,  so that detailed proofs will not be given.

We first consider the example of the Heisenberg vector $a$ inserted on the
left torus (say). The corresponding differential form is $\mathcal{F}%
_{\alpha_1,\alpha_2}^{(2)}(a;\tau_1,\tau_2,\epsilon)=Z_{\alpha_1,%
\alpha_2}^{(2)}(a,x_1\vert\tau _{1},\tau _{2},\epsilon )dx_1$. Defining $%
\nu_\alpha=\alpha_1\nu_1+\alpha_2\nu_2$, for holomorphic 1-forms $\nu_i$, we
find

\begin{theorem}
\label{theorem:Z2a_alpha} The Heisenberg vector 1-point function for a pair
of modules $M\otimes e^{\alpha_1},M\otimes e^{\alpha_2}$ is 
\begin{eqnarray}
\mathcal{F}_{\alpha_1,\alpha_2}^{(2)}(a;\tau_1,\tau_2,\epsilon) =\nu_\alpha
Z_{\alpha_1,\alpha_2}^{(2)}(\tau_1,\tau_2,\epsilon).  \label{eq: Z2a_alpha}
\end{eqnarray}
\end{theorem}
\textbf{Proof.} The proof proceeds along the same lines as Theorems \ref%
{theorem:Z2aa} and \ref{Theorem_Z2_Malpha}. We find that 
\begin{equation*}
\mathcal{F}_{\alpha_1,\alpha_2}^{(2)}(a;\tau_1,\tau_2,\epsilon)=
Z_{M}^{(1)}(\tau _{1})Z_{M}^{(1)}(\tau _{2})\sum_{D}\frac{\zeta (D)}{%
\prod_{i}e_{i}!}dx_1,
\end{equation*}
where the sum is taken over isomorphism classes of chequered configurations $%
D$ where, in this case, each configuration includes one distinguished
valence one node of type $1,x_1$. Each $D$ can be decomposed into exactly
one necklace configuration of type $\mathcal{N}_{11}^{x_1,1}$ of (\ref%
{Nabx1def}), standard configurations of the type appearing in Theorem \ref%
{Theorem_Z2_boson} and necklace contributions of type $\mathcal{N}_{ab}$ of (%
\ref{Nabdef}) as in Theorem \ref{Theorem_Z2_Malpha}. The result then follows
on applying the graphical expansion for $\nu_i(x_1)$ of (\ref{nui_graph}%
).\hfill $\square$

In a similar fashion one can generalize Theorem \ref{theorem:Z2an}
concerning the $n$-point function for $n$ Heisenberg vectors. This is
determined by the global symmetric meromorphic $n$ form given by a symmetric
(tensor) product of $\nu_\alpha$ and $\omega^{(2)}$ defined by 
\begin{equation}
\mathrm{Sym}_n\big(\omega^{(2)},\nu_\alpha\big)=
\sum_{\psi}\prod_{(r,s)}\omega^{(2)}(x_r,x_s)\prod_{(t)}\nu_\alpha(x_t),
\label{eq: bos n form alpha}
\end{equation}
where the sum is taken over the set of involutions $\psi=\ldots (rs)\ldots
(t)\ldots$ of the labels $\{1,\ldots,n\}$. Then one finds

\begin{theorem}
\label{theorem:Z2anModule} The genus two Heisenberg vector $n$-point
function for a pair of modules $M\otimes e^{\alpha_1},M\otimes e^{\alpha_2}$
is given by the global symmetric meromorphic $n$-form 
\begin{equation}
\mathcal{F}^{(2)}_{\alpha_1,\alpha_2}(a,\ldots,a;\tau_1,\tau_2,\epsilon)=%
\mathrm{Sym}_n\big(\omega^{(2)},\nu_\alpha\big) Z_{\alpha_1,\alpha_2}^{(2)}.
\label{eq: F2n alpha}
\end{equation}
\end{theorem}

Theorem \ref{theorem:Z2anModule} is a natural generalization of Corollary 4
of \cite{MT1} concerning genus one $n$-point functions for a Heisenberg
module.

Similarly to Proposition \ref{Prop:Z2omega} it follows that

\begin{proposition}
\label{Prop:Z2omega alpha} The genus two 1-point function for a pair of
modules $M\otimes e^{\alpha_1},M\otimes e^{\alpha_2}$ for the Virasoro
vector $\tilde\omega$ is 
\begin{eqnarray}
\mathcal{F}_{\alpha_1,\alpha_2}^{(2)}({\tilde\omega};\tau_1,\tau_2,%
\epsilon)= \big(\frac{1}{2}\nu_\alpha^2+ \frac{1}{12}s^{(2)} \big ) %
Z_{\alpha_1,\alpha_2}^{(2)}(\tau_1,\tau_2,\epsilon).
\label{eq: F2omega alpha}
\end{eqnarray}
\end{proposition}

Finally, let us introduce the differential operator \cite{Fa, U} 
\begin{equation}
\mathcal{D}=\frac{1}{2\pi i}\sum_{1\le i\le j\le 2}\nu_i\nu_j\frac{\partial}{%
\partial \Omega_{ij}}.  \label{eq: Del 2form}
\end{equation}
$\mathcal{D}$ maps differentiable functions on ${\mathbb{H}}_2$ to the space
of holomorphic 2-forms (spanned by $\nu_1^2, \nu_2^2, \nu_1\nu_2$) and is $%
Sp(4,\mathbb{Z})$ invariant. It follows from Theorem \ref{Theorem_Z2_Malpha}
that (\ref{eq: F2omega alpha}) can be rewritten as a Ward identity: 
\begin{equation}
\mathcal{F}_{\alpha_1,\alpha_2}^{(2)}({\tilde\omega};\tau_1,\tau_2,%
\epsilon)= Z_{M}^{(2)}(\tau_1,\tau_2,\epsilon) \big(\mathcal{D} + \frac{1}{12%
}s^{(2)} \big )e^{i\pi \alpha.\Omega.\alpha}.  \label{eq: F2omega alpha Del }
\end{equation}

% end of module discussion

\subsection{Lattice VOAs}

\label{subsect_LatticeVOA} Let $L$ be an even lattice of dimension $l$ with $%
V_{L}$ the corresponding lattice VOA. The underlying Fock space is 
\begin{equation}
V_{L}=M^{l}\otimes C[L]=\oplus _{\alpha \in L}M^{l}\otimes e^{\alpha },
\label{lattice fock space}
\end{equation}%
where $M^{l}$ is the corresponding rank $l$ Heisenberg free boson theory. We
follow Section~\ref{subsect_VOA} and \cite{MT1} concerning further notation
for lattice theories.

\medskip The general shape of $Z_{V_{L}}^{(2)}(\tau _{1},\tau _{2},\epsilon)$
is as in (\ref{Z2_def_eps}). Note that the modes of a state $u\otimes
e^{\alpha }$ map $M^{l}\otimes e^{\beta }$ to $M^{l}\otimes e^{\alpha
+\beta} $. Thus if $\alpha \neq 0$ then $Z_{V_{L}}^{(1)}(u\otimes e^{\alpha
},\tau) $ vanishes, and as a result we see that 
\begin{eqnarray}
Z_{V_{L}}^{(2)}(\tau _{1},\tau _{2},\epsilon )&=&\sum_{n\geq 0}\epsilon
^{n}\sum_{u\in M_{[n]}^{l}}Z_{V_{L}}^{(1)}(u,\tau _{1})Z_{V_{L}}^{(1)}({\bar
u},\tau _{2})  \notag \\
&=&\sum_{\alpha ,\beta \in L}\sum_{n\geq 0}\epsilon^{n}\sum_{u\in
M_{[n]}^{l}}Z_{M^{l}\otimes e^{\alpha }}^{(1)}(u,\tau _{1})Z_{M^{l}\otimes
e^{\beta }}^{(1)}({\bar u},\tau _{2}).  \label{eq: lattice part func diag}
\end{eqnarray}%
Here, $u$ ranges over any basis for $M_{[n]}^{l}$. Viewing $M^{l}\otimes
e^{\alpha }$ as a simple module for $M^{l}$ we may employ Theorem \ref%
{Theorem_Z2_Malpha} for each component to obtain

\begin{theorem}
\label{Theorem_Z2_L_eps}We have 
\begin{equation}
Z_{V_{L}}^{(2)}(\tau _{1},\tau _{2},\epsilon )=Z_{M^{l}}^{(2)}(\tau
_{1},\tau _{2},\epsilon )\theta _{L}^{(2)}(\Omega ),
\label{normalizedpartfunc}
\end{equation}%
where $\theta _{L}^{(2)}(\Omega )$ is the (genus two) Siegel theta function
associated to $L$ (e.g. \cite{Fr})%
\begin{equation}
\theta _{L}^{(2)}(\Omega )=\sum_{\alpha ,\beta \in L}\exp (\pi i((\alpha
,\alpha )\Omega _{11}+2(\alpha ,\beta )\Omega _{12}+(\beta ,\beta )\Omega
_{22})).  \label{eq: g2omega}
\end{equation}
\end{theorem}

We can similarly compute $n$-point functions for any $n$ Heisenberg vectors $%
a_1,\ldots a_l$ using Theorem \ref{theorem:Z2anModule}. We can also employ
Proposition \ref{Prop:Z2omega alpha} and the Ward identity (\ref{eq: F2omega
alpha Del }) to obtain the 1-point function for the Virasoro vector $%
\tilde\omega=\frac{1}{2}\sum_{i} a_i[-1]a_i$ as follows:

\begin{proposition}
\label{Prop: Lattice Virasoro} The Virasoro 1-point function for a lattice
VOA satisfies a genus two Ward identity 
\begin{equation}
\mathcal{F}_{V_{L}}^{(2)}({\tilde\omega};\tau_1,\tau_2,\epsilon)=
Z_{M^l}^{(2)}(\tau_1,\tau_2,\epsilon) \big(\mathcal{D} + \frac{l}{12}s^{(2)} %
\big )\theta _{L}^{(2)}(\Omega ).  \label{eq: F2omega lattice Del }
\end{equation}
\end{proposition}

The Ward identity (\ref{eq: F2omega lattice Del }) is reminiscent of some
earlier results in physics and mathematics e.g. \cite{EO, KNTY}.

We briefly discuss the automorphic properties of $Z_{V_{L}}^{(2)}(\tau
_{1},\tau _{2},\epsilon )$ and $\mathcal{F}_{L}^{(2)}({\tilde\omega}%
;\tau_1,\tau_2,\epsilon)$. There is more that one can say here, but a fuller
discussion must wait for another time \cite{MT5}. The function $\theta
_{L}^{(2)}(\Omega )$ is a Siegel modular form of weight $l/2$ (\cite{Fr})
for some subgroup of $Sp(4,\mathbb{Z})$, in particular it is holomorphic on
the Siegel upper half-space $\mathbb{H}_{2}$. From Theorems \ref%
{Theorem_period_eps}, \ref{Theorem_Z2_boson_eps_hol} and \ref%
{Theorem_Z2_L_eps}, we deduce

\begin{theorem}
\label{Theorem_Z2_lattice_eps_hol}$Z_{V_{L}}^{(2)}(\tau _{1},\tau
_{2},\epsilon )$ is holomorphic on\ the domain $\mathcal{D}^{\epsilon }$.
\hfill $\square $
\end{theorem}

We can obtain the automorphic properties of $Z_{V_{L}}^{(2)}(\tau
_{1},\tau_{2},\epsilon )$ in the same way using that for $\theta
_{L}^{(2)}(\Omega )$ together with Theorem \ref{Theorem_Z2_G}. Rather than
do this explicitly, let us introduce a variation of the partition function,
namely the \emph{normalized partition function} 
\begin{equation}
\hat{Z}_{V_{L}}^{(2)}(\tau _{1},\tau _{2},\epsilon )=\frac{%
Z_{V_{L}}^{(2)}(\tau _{1},\tau _{2},\epsilon )}{Z_{M^{l}}^{(2)}(\tau
_{1},\tau _{2},\epsilon )}.  \label{ZV_hat_eps}
\end{equation}%
Bearing in mind the convention (\ref{Omegaconvention}), what (\ref%
{normalizedpartfunc}) says is that there is a commuting diagram of
holomorphic maps 
\begin{equation}
\begin{array}{lll}
\ \hspace{0cm}\mathcal{D}^{\epsilon } & \overset{F^{\epsilon }}{%
\longrightarrow } & \hspace{0.7cm}\mathbb{H}_{2} \\ 
\hat{Z}_{V_{L}}^{(2)}\searrow &  & \swarrow \theta _{L}^{(2)} \\ 
& \hspace{0.2cm}\mathbb{C} & 
\end{array}
\label{commhol_diag}
\end{equation}

Furthermore, the $G$-actions on the two functions in question are
compatible. More precisely, if $\gamma \in G$ then we have 
\begin{eqnarray}
\hat{Z}_{V_{L}}^{(2)}(\tau _{1},\tau _{2},\epsilon )|_{l/2}\ \gamma &=&\hat{Z%
}_{V_{L}}^{(2)}(\gamma (\tau _{1},\tau _{2},\epsilon ))\det (C\Omega
+D)^{-l/2}  \notag \\
&=&\theta _{L}^{(2)}(F^{\epsilon }(\gamma (\tau _{1},\tau _{2},\epsilon
)))\det (C\Omega +D)^{-l/2}\qquad (\text{from (\ref{commhol_diag})})  \notag
\\
&=&\theta _{L}^{(2)}(\gamma (F^{\epsilon }(\tau _{1},\tau _{2},\epsilon
)))\det (C\Omega +D)^{-l/2}\qquad (\text{from Theorem \ref{TheoremGequiv}}) 
\notag \\
&=&\theta _{L}^{(2)}(\gamma \Omega )\det (C\Omega +D)^{-l/2}\qquad \ \ \ \ \
\ \ \ \ \ \ \ \ (\text{from (\ref{Omegaconvention})})  \notag \\
&=&\theta _{L}^{(2)}(\Omega )|_{l/2}\ \gamma .  \label{wtlaction}
\end{eqnarray}

For example, if the lattice $L$ is \emph{unimodular} as well as even then $%
\theta _{L}^{(2)}$ is a Siegel modular form of weight $l/2$ on the full
group $Sp(4,\mathbb{Z})$. Then (\ref{wtlaction}) informs us that 
\begin{equation*}
\hat{Z}_{V_{L}}^{(2)}(\tau _{1},\tau _{2},\epsilon )|_{l/2}\gamma =\hat{Z}%
_{V_{L}}^{(2)}(\tau _{1},\tau _{2},\epsilon ),\ \gamma \in G,
\end{equation*}%
i.e. $\hat{Z}_{V_{L}}^{(2)}(\tau _{1},\tau _{2},\epsilon )$ is automorphic
of weight $l/2$ with respect to the group $G$.

Similar remarks may be made about the \emph{normalized Virasoro 1-point
function} defined by 
\begin{equation}
\hat{\mathcal{F}}_{V_{L}}^{(2)}(\tilde\omega;\tau _{1},\tau _{2},\epsilon )=%
\frac{\mathcal{F}_{V_{L}}^{(2)}(\tilde\omega;\tau _{1},\tau _{2},\epsilon )}{%
Z_{M^{l}}^{(2)}(\tau _{1},\tau _{2},\epsilon )},  \label{Zomega_hat_eps}
\end{equation}%
which obeys the Ward identity 
\begin{equation}
\hat{\mathcal{F}}_{V_{L}}^{(2)}({\tilde\omega};\tau_1,\tau_2,\epsilon)= \big(%
\mathcal{D} + \frac{l}{12}s^{(2)} \big )\hat{Z}_{V_{L}}^{(2)}(\tau _{1},\tau
_{2},\epsilon).  \label{eq: F2omega lattice Del hat}
\end{equation}
Using the modular transformation properties of the projective connection
(e.g. \cite{Fa, U}) one finds that (\ref{eq: F2omega lattice Del hat})
enjoys the same modular properties as $\hat{Z}_{V_{L}}^{(2)}(\tau _{1},\tau
_{2},\epsilon)$ i.e.

\begin{proposition}
\label{prop: F2omega lattice mod} For the normalized Virasoro 1-point
function for a lattice VOA we have 
\begin{equation}
\hat{\mathcal{F}}_{V_{L}}^{(2)}({\tilde\omega};\tau_1,\tau_2,%
\epsilon)|_{l/2}\gamma= \big(\mathcal{D} + \frac{l}{12}s^{(2)} \big )\big(%
\hat{Z}_{V_{L}}^{(2)}(\tau _{1},\tau _{2},\epsilon )|_{l/2}\gamma\big),
\label{eq: F2omega lattice Del mod}
\end{equation}
for $\gamma \in G$.
\end{proposition}

\subsection{Rank Two Fermion Vertex Super Algebra and the Genus Two Riemann
Theta Series}

As a last application of Theorem \ref{Theorem_Z2_Malpha}, we briefly
consider the rank two fermion Vertex Operator Super Algebra (VOSA) $V=V(H,{}+%
\frac{1}{2})^2$. $V$ can be decomposed in terms of a Heisenberg subVOA
generated by a Heisenberg state $a$ and irreducible modules $M\otimes e^m$
for $m\in$ e.g. \cite{Ka}. One can construct orbifold $n$-point functions
for a pair $g,h$ of commuting $V$ automorphisms generated by $a(0)$ \cite%
{MTZ}. In particular, consider the 1-point function (which is non-vanishing
only for $u\in M$) for a $g$-twisted sector for $g=e^{-2\pi i\lambda a(0)}$
together with an automorphism $h=e^{2\pi i\mu a(0)}$ (for real $\lambda,\mu$%
) which can be expressed as (op. cite.) 
\begin{eqnarray}
Z^{(1)}_V((g,h);u,\tau)&=&\mathrm{Tr}_V(h o(u)q^{L(0)+\lambda^2/2+\lambda
a(0)-1/24})  \notag \\
&=& \sum_{m\in \mathbb{Z}}e^{2\pi im\mu }\mathrm{Tr}_{M\otimes
e^{m+\lambda}}(o(u)q^{L(0)-1/24}),  \label{eq: Zgh part}
\end{eqnarray}
utilizing the Heisenberg decomposition. In particular, the orbifold
partition function is expressed in terms of the Jacobi theta series 
\begin{eqnarray}
Z^{(1)}_V((g,h);\tau) &=&\frac{e^{-2\pi i\lambda\mu}}{\eta (\tau )}\vartheta %
\left[ 
\begin{array}{c}
\lambda \\ 
\mu%
\end{array}%
\right] (\tau ),  \notag \\
\vartheta \left[ 
\begin{array}{c}
\lambda \\ 
\mu%
\end{array}%
\right] (\tau )&=& \sum_{m\in \mathbb{Z}}e^{i\pi (m+\lambda)^2\tau+2\pi i
(m+\lambda)\mu }.  \label{Z_lattice}
\end{eqnarray}

Similarly to (\ref{eq: alpha part func}), it is natural to define the genus
two orbifold partition function for a pair of $g_i$-twisted sectors together
with commuting automorphisms $h_i$ parameterized by $\lambda_i,\mu_i$ for $%
i=1,2$ with 
\begin{equation}
Z_V^{(2)}((g_i,h_i);\tau _{1},\tau _{2},\epsilon )=\sum_{n\geq 0}\epsilon
^{n}\sum_{u\in M_{[n]}}Z^{(1)}_V((g_1,h_1);u,\tau _{1}) Z^{(1)}_V((g_2,h_2);{%
\bar u},\tau _{2}),  \label{eq: gh_i part func}
\end{equation}
where $u$ ranges over any basis for $M_{[n]}$. A more detailed description
of this and an alternative fermionic VOSA approach to this will be described
elsewhere \cite{TZ}. Here we decompose the genus one 1-point functions of (%
\ref{eq: gh_i part func}) in terms of Heisenberg modules $M\otimes
e^{m_i+\lambda_i}$ to find, in the notation of (\ref{eq: alpha part func}),
that 
\begin{equation}
Z_V^{(2)}((g_i,h_i);\tau _{1},\tau _{2},\epsilon )=\sum_{m\in \mathbb{Z}^2}
e^{2\pi i m.\mu }Z_{m_1+\lambda_1,m_2+\lambda_2}^{(2)}(\tau _{1},\tau
_{2},\epsilon ),  \label{eq: Heis expansion}
\end{equation}
where here $\lambda=(\lambda_1,\lambda_2),\mu=(\mu_1,\mu_2)\in \mathbb{R}^2$
and $m=(m_1,m_2)\in \mathbb{Z}^2$. Theorem \ref{Theorem_Z2_Malpha} implies

\begin{theorem}
\label{Theorem_Z2_orbifold} We have 
\begin{equation}
Z_V^{(2)}((g_i,h_i);\tau _{1},\tau_{2},\epsilon )=e^{-2\pi
i\lambda.\mu}Z_{M}^{(2)}(\tau _{1},\tau _{2},\epsilon )\theta^{(2)}\left[ 
\begin{array}{l}
\lambda \\ 
\mu%
\end{array}
\right] (\Omega ),  \label{Z2_orbifoldpartfunc}
\end{equation}%
for genus two Riemann theta function (e.g. \cite{Mu})%
\begin{equation}
\theta^{(2)}\left[ 
\begin{array}{l}
\lambda \\ 
\mu%
\end{array}
\right] {(\Omega )} =\sum_{m\in \mathbb{Z}^2}
e^{i\pi(m+\lambda).\Omega.(m+\lambda)+2\pi i(m+\lambda).\mu}.
\end{equation}
\end{theorem}

As already described for lattice VOAs, one can similarly obtain a Ward
identity for the Virasoro 1-point function analogous to (\ref{eq: F2omega
lattice Del }) and (\ref{eq: F2omega lattice Del hat}) and analyze the
modular properties of (\ref{Z2_orbifoldpartfunc}) and the Virasoro 1-point
function under the action of $G$.

\section{Appendix - A Product Formula}

\label{sect_Appendix} Here we continue the discussion initiated in
Subsection~\ref{subsect_cycles}, with a view to proving Proposition \ref%
{Prop_Om12_R21expansion}. Consider a set of independent (non-commuting)
variables $x_{i}$ indexed by the elements of a finite set $I=\{1,\ldots ,N\}$%
. The set of all distinct monomials $x_{i_{1}}\ldots x_{i_{n}}(n\geq 0)$ may
be considered as a basis for the tensor algebra associated with an $N$
dimensional vector space. Call $n$ the degree of the monomial $%
x_{i_{1}}\ldots x_{i_{n}}$.

Let $\rho =\rho _{n}$ be the standard cyclic permutation which acts on
monomials of degree $n$ via $\rho :x_{i_{1}}\ldots x_{i_{n}}\mapsto
x_{i_{n}}x_{i_{1}}\ldots x_{i_{n-1}}$. The \emph{rotation group} of a given
monomial $x=x_{i_{1}}\ldots x_{i_{n}}$ is the subgroup of $\langle \rho
_{n}\rangle $ that leaves $x$ invariant. Call $x$ \emph{rotationless} in
case its rotation group is trivial. Let us say that two monomials $x,y$ of
degree $n$ are \emph{equivalent} in case $y=\rho _{n}^{r}(x)$ for some $r\in
Z$, and denote the corresponding equivalence class by $(x)$. We call these 
\emph{cycles}. Note that equivalent monomials have the same rotation group,
so we may meaningfully refer to the rotation group of a cycle. In
particular, a \emph{rotationless cycle} is a cycle whose representative
monomials are themselves rotationless. Let $C_{n}$ be the set of
inequivalent cycles of degree $n$.

It is convenient to identify a cycle $(x_{i_{1}}\ldots x_{i_{n}})$ with a 
\emph{cyclic labeled graph} or \emph{labeled polygon}, that is, a graph with 
$n$ vertices labeled $x_{i_{1}},\ldots ,x_{i_{n}}$ and with edges $%
x_{i_{1}}x_{i_{2}},\ldots ,x_{i_{n-1}}x_{i_{n}},x_{i_{n}}x_{i_{1}}$. We will
sometimes afflict the graph with one of the two canonical orientations.

A cycle is rotationless precisely when its graph admits no non-trivial
rotations (a rotation now being an orientation-preserving automorphism of
the graph which preserves labels of nodes).

Let $\mathcal{M}(I)$ be the (multiplicative semigroup generated by) the
rotationless cycles in the symbols $x_{i},i\in I$. There is an injection 
\begin{equation}
\iota :\bigcup_{n\geq 0}C_{n}\longrightarrow \mathcal{M}(I)
\label{eq: cycleinjection}
\end{equation}%
defined as follows. If $(x)\in C_{n}$ has rotation group of order $r$ then $%
r|n$ and there is a rotationless monomial $y$ such that $x=y^{r}$. We then
map $(x)\mapsto (y)^{r}$. It is readily verified that this is well-defined.
In this way, each cycle is mapped to a power of a rotationless cycle in $%
\mathcal{M}(I)$. A typical element of $\mathcal{M}(I)$ is uniquely
expressible in the form 
\begin{equation}
p_{1}^{f_{1}}p_{2}^{f_{2}}\ldots p_{k}^{f_{k}}  \label{eq:RF0}
\end{equation}%
where $p_{1},\ldots ,p_{k}$ are distinct rotationless cycles and $%
f_{1},\ldots ,f_{k} $ are non-negative integers. We call (\ref{eq:RF0}) the 
\textit{reduced form} of an element in $\mathcal{M}(I)$. A general element
of $\mathcal{M}(I)$ is then essentially a labeled graph, each of whose
connected components are rotationless labeled polygons as discussed in
Subsection~\ref{subsect_cycles}.

Now consider a second finite set $T$ together with a map 
\begin{equation}
F:T\longrightarrow I.  \label{eq:RF1}
\end{equation}%
Thus elements of $I$ label elements of $T$ via the map $F$. $F$ induces a
natural map 
\begin{equation*}
\overline{F}:\Sigma (T)\longrightarrow \mathcal{M}(I)
\end{equation*}%
from the symmetric group $\Sigma (T)$ as follows. For an element $\tau \in
\Sigma (T)$, write $\tau $ as a product of disjoint cycles $\tau =\sigma
_{1}.\sigma _{2}\ldots $. We set $\overline{F}(\tau )=\overline{F}(\sigma
_{1})\overline{F}(\sigma _{2})\ldots $, so it suffices to define $\overline{F%
}(\sigma )$ for a cycle $\sigma =(s_{1}s_{2}\ldots )$ with $%
s_{1},s_{2},\ldots \in T$. In this case we set 
\begin{equation*}
\overline{F}(\sigma )=\iota ((x_{F(s_{1})}x_{F(s_{2})}\ldots ))
\end{equation*}%
where $\iota $ is as in (\ref{eq: cycleinjection}). When written in the form
(\ref{eq:RF0}), we call $\overline{F}(\tau )$ the \textit{reduced $F$-form}
of $\tau $.

For $i\in {I}$, let $s_{i}=|F^{-1}(i)|$ be the number of elements in $T$
with label $i$. So the number of elements in $T$ is equal to ${\sum }_{i\in
I}s_{i}$. We say that two elements $\tau _{1},\tau _{2}\in \Sigma (T)$ are $%
F $-equivalent if they have the same reduced $F$-form, i.e. $\overline{F}%
(\tau _{1})=\overline{F}(\tau _{2})$. We will show that each equivalence
class contains the same number of elements. Precisely,

\begin{lemma}
\label{Lemma_Fequiv} Each $F$-equivalence class contains precisely ${\prod }%
_{i\in I}s_{i}$ elements. In particular, the number of $F$-equivalence
classes is $|T|!/{\prod }_{i\in I}s_{i}$.
\end{lemma}

\noindent \textbf{Proof.} An element $\tau \in \Sigma (T)$ may be
represented uniquely as 
\begin{equation*}
\left( 
\begin{array}{cccc}
0 & 1 & \cdots & M \\ 
\tau (0) & \tau (1) & \cdots & \tau (M)%
\end{array}%
\right)
\end{equation*}%
so that 
\begin{equation*}
\overline{F}(\tau )=\left( 
\begin{array}{cccc}
F(0) & F(1) & \cdots & F(M) \\ 
F(\tau (0)) & F(\tau (1)) & \cdots & F(\tau (M))%
\end{array}%
\right)
\end{equation*}%
with an obvious notation. Exactly $s_{i}$ of the $\tau (j)$ satisfy 
\begin{equation*}
\overline{F}(\tau (j))=x_{i}
\end{equation*}%
so that there are $\prod_{i\in I}s_{i}$ choices of $\tau $ which have a
given image under $\overline{F}$. The Lemma follows. \hfill ${\square }$

\medskip The next results employs notation introduced in Subsections~\ref%
{subsect_cycles} and \ref{subsect_Omega_graph}.

\begin{lemma}
\label{Lemma_L21_expansion} We have 
\begin{equation}
(I-M_{1}M_{2})^{-1}(1,1)=(1-\sum_{L\in \mathcal{L}_{21}}\zeta (L))^{-1}.
\label{eq: matrixsum}
\end{equation}
\end{lemma}

\noindent As before, the left-hand-side of (\ref{eq: matrixsum}) means $%
\sum_{n\geq 0}(M_{1}M_{2})^{n}(1,1)$. It is a certain power series with
entries being quasi-modular forms.

\medskip \noindent \textbf{Proof of Lemma.} We have 
\begin{equation}
(M_{1}M_{2})^{n}(1,1)=\sum M_{1}(1,k_{1})M_{2}(k_{1},k_{2})\ldots
M_{2}(k_{2n-1},1)  \label{eq: matrixentry}
\end{equation}%
where the sum ranges over all choices of positive integers $k_{1},\ldots
,k_{2n-1}$. Such a choice corresponds to a (isomorphism class of) chequered
cycle $L$ with $2n$ nodes and with at least one distinguished node, so that
the left-hand-side of (\ref{eq: matrixsum}) is equal to 
\begin{equation*}
\sum_{L}\zeta (L)
\end{equation*}%
summed over all such $L$. We can formally write $L$ as a product $%
L=L_{1}L_{2}\ldots L_{p}$ where each $L_{i}\in \mathcal{L}_{21}$. This
indicates that $L$ has $p$ distinguished nodes and that the $L_{i}$ are the
edges of $L $ between consecutive distinguished nodes, which can be
naturally thought of as chequered cycles in $\mathcal{L}_{21}$. Note that in
the representation of $L$ as such a product, the $L_{i}$ do not commute
unless they are equal, moreover $\zeta $ is multiplicative. Then 
\begin{equation*}
(I-M_{1}M_{2})^{-1}(1,1)=\sum_{L_{i}\in \mathcal{L}_{21}}\zeta (L_{1}\ldots
L_{p})=(1-\sum_{L\in \mathcal{L}_{21}}\zeta (L))^{-1}
\end{equation*}%
as required. \hfill$\square $

\begin{proposition}
\label{Prop_R21expansion}We have 
\begin{equation}
(I-M_{1}M_{2})^{-1}(1,1)=\prod_{L\in \mathcal{R}_{21}}(1-\zeta (L))^{-1}
\label{eq: matrixprod}
\end{equation}
\end{proposition}

\noindent \textbf{Proof.} By Lemma \ref{Lemma_L21_expansion} we have 
\begin{equation}
(I-M_{1}M_{2})^{-1}(1,1)=\sum m(e_{1},\ldots ,e_{k})\zeta
(L_{1})^{e_{1}}\ldots \zeta (L_{k})^{e_{k}}  \label{eq: mult}
\end{equation}%
where the sum ranges over distinct elements $L_{1},\ldots L_{k}$ of $%
\mathcal{L}_{21}$ and all $k$-tuples of non-negative integers $e_{1},\ldots
,e_{k}$, and where the multiplicity is 
\begin{equation*}
m(e_{1},\ldots ,e_{k})=\frac{(\sum_{i}e_{i})!}{\prod_{i}(e_{i}!)}.
\end{equation*}%
Let $S$ be the set consisting of $e_{i}$ copies of $L_{i},1\leq i\leq k$,
let $I$ be the integers between $1$ and $k$, and let $F:S\longrightarrow I$
be the obvious labelling map. A reduced $F$-form is then an element of $%
\mathcal{M}(I)$ where the variables $x_{i}$ are now the $L_{i}$. The free
generators of $\mathcal{M}(I)$, i.e. rotationless cycles in the $x_{i}$, are
naturally identified \emph{precisely} with the elements of $\mathcal{R}_{21}$%
, and Lemma \ref{Lemma_Fequiv} implies that each element of $\mathcal{M}(I)$
corresponds to just one term under the summation in (\ref{eq: mult}). Eqn.(%
\ref{eq: matrixprod}) follows immediately from this and the multiplicativity
of $\zeta $, and the Proposition is proved. \hfill $\square $

\end{document}